\documentstyle[amssymb]{amsart}

\newif\ifTwelve

\ifTwelve \let\fixedline\newline \else \let\fixedline\relax \fi

\ifx\mathrm\undefined\let\mathrm\bf\fi
\ifx\mathbf\undefined\let\mathbf\bf\fi

\let\leq\leqslant \let\geq\geqslant
\let\tsize\textstyle \def\Sum{\sum\limits}

\let\Box\square
\let\al\alpha
\let\bt\beta
\let\gm\gamma \let\Gm\Gamma
 
 \let\eps\varepsilon \let\epsilon\eps

\let\la\lambda \let\La\Lambda
\let\om\omega \let\Om\Omega
 \let\phi\varphi

\newcommand{\half}{\frac12}
\newcommand{\Z}{{\Bbb Z}}

\newcommand{\R}{{\Bbb R}}
\newcommand{\C}{{\Bbb C}}

\newcommand{\Ref}[1]{{$($\ref{#1}$)$}}
\newcommand{\bean}{\begin{eqnarray}}
\newcommand{\eean}{\end{eqnarray}}
\newcommand{\be}{\begin{displaymath}}
\newcommand{\ee}{\end{displaymath}}
\newcommand{\bea}{\begin{eqnarray*}}
\newcommand{\eea}{\end{eqnarray*}}

\newcommand{\h}{{{\frak h\,}}}
\newcommand{\Id}{{\operatorname{Id}}}

\newcommand{\noi}{\noindent}
\newcommand{\vs}{\vspace{.5\baselineskip}}
\newcommand{\Mu}{{\mathrm M}}
\newenvironment{proof}{\noindent{\it Proof\/}:}{$\;\Box$\par\vs}
\newenvironment{definition}
{\noindent{\bf Definition\/}:}{\par\vs}
\newenvironment{example}
{\noindent{\bf Example\/}:}{\par\vs}
\newtheorem%
{thm}{Theorem}
\newtheorem%
{proposition}[thm]{Proposition}
\newtheorem%
{lemma}[thm]{Lemma}
\newtheorem%
{lemmadef}[thm]{Lemma-Definition}
\newtheorem%
{corollary}[thm]{Corollary}
\newtheorem%
{conjecture}[thm]{Conjecture}

\newcommand{\End}{{\operatorname{End\,}}}

\newcommand{\Sym}{{\operatorname{Sym\,}}}
\newcommand{\Fun}{{\operatorname{Fun\,}}}
\renewcommand{\Im}{{\operatorname{Im\,}}}

\def\nli#1{\vsk#1>\noindent\ignorespaces}
\def\nl{\nli0}
\def\vsk#1>{\vskip#1\baselineskip}
\def\ftext#1{{\let\thefootnote\relax\footnotetext{\vsk-.8>\noindent #1}}}
\newenvironment{abst}{\begingroup\small\narrower\noindent{\bf Abstract.}
\enspace\ignorespaces}{\endgraf\endgroup}

\title[Monodromy of qKZB difference equations]
{Monodromy of solutions of the elliptic quantum\\[4pt]
Knizhnik-Zamolodchikov-Bernard\\[4pt]
difference equations }

\author[G. Felder, V. Tarasov and A. Varchenko]
{G.\ Felder$^{\,\star}$, V.\ Tarasov$^{\,*}$ and
A.\ Varchenko$^{\,\diamond}$}

\begin{document}
\maketitle

\begin{center}
{\it
$^\star$Department of Mathematics, ETH Zentrum, CH-8092 Z\"urich, Switzerland
\vsk.5>
$^*$Department of Mathematics, Faculty of Science, Osaka University
\vsk0>
Toyonaka, Osaka 560, Japan
\vsk.5>
$^\diamond$Department of Mathematics, University of North Carolina
at Chapel Hill
\vsk0>
Chapel Hill, NC 27599-3250, USA}
\end{center}

\ftext{\small{\normalsize\sl $^\star$Email\/{\rm:} felder@@math.ethz.ch}\nl
\hphantom{$*$}Supported in part by NSF grant DMS-9400841\nli{.2}
{\normalsize\sl $^*$Email\/{\rm:} vt@@math.sci.osaka-u.ac.jp}\nl
\hphantom{$*$}On leave of absense from St.\,Petersburg Branch of
Steklov Mathematical Institute\nli{.2}
{\normalsize\sl $^\diamond$Email\/{\rm:} av@@math.unc.edu}\nl
\hphantom{$*$}Supported in part by NSF grant DMS-9501290}

\vsk1.5>
\centerline{April, 1997}
\vsk1.8>

\begin{abst}
The elliptic quantum Knizhnik--Zamolodchikov--Bernard (qKZB) difference
equations associated to the elliptic quantum group $E_{\tau,\eta}(sl_2)$
is a system of difference equations with values in a tensor product of
representations of the quantum group and defined in terms of the elliptic
$R$-matrices associated with pairs of representations of the quantum group.
In this paper we solve the qKZB equations in terms of elliptic hypergeometric
functions and decribe the monodromy properties of solutions. It turns out that
the monodromy transformations of solutions are described in terms of elliptic
$R$-matrices associated with pairs of representations of the "dual" elliptic
quantum group $E_{p,\eta}(sl_2)$, where $p$ is the step of the difference
equations. Our description of the monodromy is analogous to the Kohno-Drinfeld
description the monodromy group of solutions of the KZ differential equations
associated to a simple Lie algebra in terms of the corresponding quantum group.
\end{abst}
\vsk>
\vsk0>

\thispagestyle{empty}

\section{Introduction}

In this paper we solve the system of elliptic quantum
Knizhnik--Zamolodchikov--Bernard (qKZB) difference equations associated with
the elliptic quantum group $E_{\tau,\eta}(sl_2)$ and describe the monodromy
properties of solutions.

The qKZB equations \cite{F} are a quantum deformation of the KZB differential
equations obeyed by correlation functions of the Wess--Zumino--Witten model on
tori. The qKZB equations have the form
\be
\Psi(z_1,\dots,z_j+p,\dots,z_n)=
K_j(z_1,\dots,z_n;\tau,\eta,p)\Psi(z_1,\dots,z_n).
\ee
The unknown function $\Psi$ takes values in a space of vector valued functions
of a complex variable $\lambda$, and the $K_j$ are difference operators in
$\lambda$. The parameters of this system of equations are $\tau$ (the period
of the elliptic curve), $\eta$ (``Planck's constant''), $p$ (the step) and $n$
``highest weights'' $\Lambda_1,\dots,\Lambda_n\in\C$. The operators $K_j$ are
expressed in terms of $R$-matrices of the elliptic quantum group
$E_{\tau,\eta}(sl_2)$.

In the trigonometric limit $\tau\to i\infty$, the qKZB equations reduce to
the trigonometric qKZ equations \cite{FR} obeyed by correlation functions of
statistical models and form factors of integrable quantum field theories in
1+1 dimensions \cite{JM, S}.

The KZB equations can be obtained in the semiclassical limit: $\eta\to 0$,
$p\to 0$, $p/\eta$ finite.

When the step $p$ of the qKZB equations goes to zero (with the other parameters
fixed) our construction gives common eigenfunctions of the $n$ commuting
operators $K_j(z_1,\dots,z_n;\tau,\eta,0)$ in the form of the Bethe ansatz
\cite{FTV}. These difference operators are closely related to the transfer
matrices of IRF models of statistical mechanics \cite{F, FV2}.

Our first main result is a construction of solutions of the qKZB equations in
the form of multidimensional elliptic hypergeometric intergals.

Our second main result is a description of the monodromy properties of
solutions. We show that the qKZB equations can be considered as difference
equations on the product of several copies of the elliptic curve with modulus
$\tau$ with values in a suitable vector bundle. Therefore, a natural question
is to describe the monodromy of solutions with respect to shifts of arguments
of solutions by periods of the elliptic curve. It turns out that the monodromy
transformations of solutions are described in terms of $R$-matrices associated
with pairs of representations of the "dual" elliptic quantum group
$E_{p,\eta}(sl_2)$, where $p$ is the step of the difference equations.

Our description of the monodromy is analogous to the Kohno-Drinfeld description
\cite{K, D} of the monodromy group of solutions of the KZ differential
equations associated to a simple Lie algebra in terms of the corresponding
quantum group.

The results of this paper are parallel to the results on solutions of the
rational and trigonometric qKZ equations of \cite{V, TV1, TV2}, which are based
on the representation theory of the Yangian $Y(sl_2)$ and the affine quantum
universal enveloping algebra $U_q(\widehat{sl_2})$, respectively.
In particular, in \cite{TV1} the monodromy of the qKZ difference equations
associated with the Yangian $Y(sl_2)$ is descibed in terms of the affine
quantum group $U_q(\widehat{sl_2})$, where the parameter $q$ is connected with
the step $p$ of the equations by $q=e^{\pi i/p}$. In \cite{TV2} the monodromy
of the qKZ difference equations associated with the affine quantum group
$U_q(\widehat{sl_2})$ is descibed in terms of the elliptic quantum group
$E_{\tau,\eta}(sl_2)$ where the parameters $\tau$ and $\eta$ of the elliptic
quantum group are connected with the parameters $q$ and the multiplicative step
$p$ of the equations by relations $p=e^{2 \pi i \tau}$ and $q=e^{-2\pi i\eta}$.

The paper is organized as follows. We begin by introducing the notion of
$R$-matrices and the qKZB equations in Section \ref{sqkzb}. The geometric
construction of $R$-matrices is given in Section \ref{ssafs}. At the end of
that section we show how to obtain representations of $E_{\tau,\eta}(sl_2)$
in this way, and give some explicit formulae for $R$-matrix elements.

In Section \ref{transf} we describe transformation properties of the qKZB
equations with respect to shifts of the arguments by $\tau$ and $1$ and show
that the qKZB equations can be considered as difference equations on a product
of several copies of the elliptic curve with modulus $\tau$.

In Section \ref{ssqkzb} we describe formal integral representations of
solutions of the qKZB equations. We constuct solutions of the qKZB equations
in Section \ref{integration}. The monodromy properties of solutions are
described in Section \ref{Monodromy}.

\section{$R$-matrices, qKZB equations }\label{sqkzb}

\subsection{$R$-matrices}
The qKZB equations are given in terms of $R$-matrices of elliptic quantum
groups. In the $sl_2$ case, these $R$-matrices have the following properties.
Let $\h=\C h$ be a one-dimensional Lie algebra with generator $h$. For each
$\Lambda\in\C$ consider the $\h$-module $V_\Lambda=\oplus_{j=0}^\infty\C e_j$,
with $he_j=(\Lambda-2j)e_j$. For each pair $\Lambda_1$, $\Lambda_2$ of complex
numbers we have a meromorphic function, called the $R$-matrix,
$R_{\Lambda_1,\Lambda_2}(z,\lambda)$ of two complex variables, with values in
$\End(V_{\Lambda_1}\otimes V_{\Lambda_2})$.

The main properties of the $R$-matrices are
\begin{enumerate}
\item[I.] The zero weight property: for any $\Lambda_i$, $z,\lambda$,
$[R_{\Lambda_1,\Lambda_2}(z,\lambda),h^{(1)}+h^{(2)}]=0$.
\item[II.] For any $\Lambda_1,\Lambda_2,\Lambda_3$, the dynamical Yang--Baxter
equation
\bea
R_{\Lambda_1,\Lambda_2}(z,\lambda-2\eta h^{(3)})^{(12)}
R_{\Lambda_1,\Lambda_3}(z+w,\lambda)^{(13)}
R_{\Lambda_2,\Lambda_3}(w,\lambda-2\eta h^{(1)})^{(23)}
\\
{}=
R_{\Lambda_2,\Lambda_3}(w,\lambda)^{(23)}
R_{\Lambda_1,\Lambda_3}(z+w,\lambda-2\eta h^{(2)})^{(13)}
R_{\Lambda_1,\Lambda_2}(z,\lambda)^{(12)},
\eea
holds in $\End(V_{\Lambda_1}\otimes V_{\Lambda_2}\otimes V_{\Lambda_3})$
for all $z,w,\lambda$.
\item[III.] For all $\Lambda_1$, $\Lambda_2$, $z,\lambda$,
$R_{\Lambda_1,\Lambda_2}(z,\lambda)^{(12)}
R_{\Lambda_2,\Lambda_1}(-z,\lambda)^{(21)}=\Id$.
This property is called the ``unitarity''.
\end{enumerate}

We use the following notation: if $X\in\End(V_i)$, then we denote by
$X^{(i)}\in\End(V_1\otimes\dots\otimes V_n)$ the operator
$\cdots\otimes\Id\otimes X\otimes\Id\otimes\cdots$, acting non-trivially on
the $i$th factor of a tensor product of vector spaces, and if
$X=\sum X_k\otimes Y_k\in\End(V_i\otimes V_j)$, then we set
$X^{(ij)}=\sum X_k^{(i)}Y_k^{(j)}$. If $X(\mu_1,\dots,\mu_n)$ is a function
with values in $\End(V_1\otimes\dots\otimes V_n)$, then
$X(h^{(1)},\dots,h^{(n)})v=X(\mu_1,\dots,\mu_n)v$ if $h^{(i)}v=\mu_iv$,
for all $i=1,\dots,n$.

For each $\tau$ in the upper half plane and generic $\eta\in\C$ (``Planck's
constant'') a system of $R$-matrices $R_{\Lambda_1,\Lambda_2}(z,\lambda)$
obeying {I\,--\,III} was constructed in \cite{FV1}. They are characterized by
an intertwining property with respect to the action of the elliptic quantum
group $E_{\tau,\eta}(sl_2)$ on tensor products of evaluation Verma modules.

\subsection{qKZB equations}
Fix the parameters $\tau,\eta$. Fix also $n$ complex numbers
$\Lambda_1,\dots,\Lambda_n$ and an additional parameter $p\in\C$.
Let $V=V_{\Lambda_1}\otimes\cdots\otimes V_{\Lambda_n}$. The kernel of
$h^{(1)}+\dots+h^{(n)}$ on $V$ is called the zero-weight space and is denoted
$V[0]$. More generally, we write $V[\mu]$ for the eigenspace of $\sum h^{(i)}$
with eigenvalue $\mu$. The qKZB equations are difference equations for
a function $\Psi(z_1,\dots,z_n,\lambda)$ of $n$ complex variables
$z_1,\dots,z_n$ with values in the space of meromorphic functions
$\Fun(V[0])$ of a complex variable $\lambda$ with values in $V[0]$.

The qKZB equations \cite{F} have the form
\bean\label{qKZB}
\Psi(z_1,\dots,z_j+p,\dots,z_n)=
R_{j,j-1}(z_j\!-\!z_{j-1}+p)\cdots R_{j,1}(z_j\!-\!z_{1}+p)
\\[2pt]
\Gamma_j
R_{j,n}(z_j\!-\!z_n)\cdots,R_{j,j+1}(z_j\!-\!z_{j+1})
\Psi(z_1,\dots,z_n)\notag
\eean
Here $R_{k,l}(z)$ is the operator of multiplication by
\be
\tsize R_{\Lambda_k,\Lambda_l}(z,\lambda-2\eta\Sum_{\tsize{j=1\atop j\neq k}}
^{l-1}h^{(j)})^{(k,l)}
\ee
acting on the $k$th and $l$th factor of the tensor product,
and $\Gamma_j$ is the linear difference operator such that
$\Gamma_j\Psi(\lambda)=\Psi(\lambda-2\eta\mu)$ if $h^{(j)}\Psi= \mu\Psi$.

The consistency of these equations follows from I--III. In other words,
the qKZB equations may be viewed as the equation of horizontality for a flat
discrete connection on a trivial vector bundle with fiber $\Fun(V[0])$ over
an open subset of $\C^n$.

\subsection{Finite-dimensional representations}
If $\Lambda$ is a nonnegative integer, $V_\Lambda$ contains the subspace
$SV_\Lambda=\oplus_{j=\Lambda+1}^\infty\C e_j$ with the property that, for any
$\Mu$, $SV_\Lambda\otimes V_\Mu$ and $V_\Mu\otimes SV_\Lambda$ are preserved by
the $R$-matrices $R_{\Lambda,\Mu}(z,\lambda)$ and $R_{\Mu,\Lambda}(z,\lambda)$,
respectively, see \cite{FV1} and Theorem \ref{tfd}. Let
$L_\Lambda=V_\Lambda/SV_\Lambda$, $\Lambda\in\Z_{\geq 0}$. Then, in particular,
for any nonnegative integers $\Lambda$ and $\Mu$, $R_{\Lambda,\Mu}(z,\lambda)$
induces a map, also denoted by $R_{\Lambda,\Mu}(z,\lambda)$, on the
finite-dimensional space $L_\Lambda\otimes L_\Mu$.

The simplest nontrivial case is $\Lambda=\Mu=1$. Then $R_{1,1}(z,\lambda)$
is defined on a four-dimensional vector space and coincides with
the {\em fundamental} $R$-matrix, the matrix of structure constants of
the elliptic quantum group $E_{\tau,\eta}(sl_2)$.

In any case, if $\Lambda_1,\dots,\Lambda_n$ are nonnegative integers, we can
consider the qKZB equations \Ref{qKZB} on functions with values in the zero
weight space of $L_{\Lambda_1}\otimes \cdots\otimes L_{\Lambda_n}$.

The results of this paper obtained for the solutions with values in
$\otimes_j V_{\Lambda_j}$ extend to this case: let
$\pi:\otimes_{j=1}^nV_{\Lambda_j}\to\otimes_{j=1}^{n} L_{\Lambda_j}$ denote
the canonical projection.

\begin{lemma}
Let $\Psi(z_1,\dots,z_n)$ be a solution of the qKZB equations with values
in $V[0]=V_{\Lambda_1}\otimes\cdots\otimes V_{\Lambda_n}[0]$.
Then $\pi\circ\Psi(z_1,\dots,z_n)$ is a solution of the qKZB equations with
values in $L[0]=L_{\Lambda_1}\otimes\cdots\otimes L_{\Lambda_n}[0]$.
\end{lemma}


\section{Modules over the elliptic quantum group as function
spaces}\label{ssafs}

In this section we realize the spaces dual to tensor products of evaluation
Verma modules over $E_{\tau,\eta}(sl_2)$ as spaces of functions.
The $R$-matrices are then constructed geometrically.

Let us fix complex parameters $\tau$, $\eta$ with $\Im\tau>0$, and complex
numbers $\Lambda_1,\dots, \Lambda_n$. We set $a_i=\eta\Lambda_i$,
$i=1,\dots,n$.

\subsection{A space of symmetric functions}\label{sssym}
Introduce a space of functions with an action of the symmetric group.
Recall that the Jacobi theta function
\begin{equation}
\theta(t)=-\sum_{j\in\Z}
e^{\pi i(j+\half)^2\tau+2\pi i(j+\half)(t+\half)},
\notag
\end{equation}
has multipliers $-1$ and $-\exp(-2\pi it-\pi i\tau)$ as $t\to t+1$ and
$t\to t+\tau$, respectively. It is an odd entire function whose zeros are
simple and lie on the lattice $\Z+\tau \Z$. It has the product formula
\be
\theta(t)\,=\,2e^{\pi i\tau/4}\sin(\pi t)
\prod_{j=1}^\infty(1-q^j)(1-q^je^{2\pi it})(1-q^je^{-2\pi it}),
\qquad q=e^{2\pi i\tau}.
\ee

\begin{definition}
For complex numbers $a_1,\dots,a_n$, $z_1,\dots,z_n$, $\lambda$,
let $\tilde F^m_{a_1,\dots,a_n}(z_1,\dots,z_n,\lambda)$ be the space of
meromorphic functions $f(t_1,\dots,t_m)$ of $m$ complex variables such that
\begin{enumerate}
\item[(i)] $\prod_{i<j}\theta(t_i-t_j+2\eta)
\prod_{i=1}^m\prod_{k=1}^n\theta(t_i-z_k-a_k)f$
is a holomorphic function on $\C^m$.
\item[(ii)] $f$ is periodic with period 1 in each of its arguments and
\be
f(\cdots,t_j+\tau,\cdots)
=e^{-2\pi i(\lambda+4\eta j-2\eta)}f(\cdots,t_j,\cdots),
\ee
for all $j=1,\dots,m$.
\end{enumerate}
\end{definition}
\noi
The symmetric group $S_m$ acts on
$\tilde F^m_{a_1,\dots,a_n}(z_1,\dots,z_n,\lambda)$ so that the transposition
of $j$ and $j+1$ acts as
\be
s_jf(t_1,\dots,t_m)=f(t_1,\dots,t_{j+1},t_j,\dots,t_m)
\frac{\theta(t_{j}-t_{j+1}-2\eta)}
{\theta(t_{j}-t_{j+1}+2\eta)}\,.
\ee

\begin{definition}
For any $m\in\Z_{>0}$, let $F^m_{a_1,\dots,a_n}(z_1,\dots,z_n,\lambda)=
\tilde F^m_{a_1,\dots,a_n} (z_1,\dots,z_n,\lambda)^{S_m}$ be the space of
$S_m$-invariant functions. If $m=0$, then we set
$F^0_{a_1,\dots,a_n}(z_1,\dots,z_n,\lambda)=\C$. We denote by $\Sym$
the symmetrization operator $\Sym=\sum_{s\in S_m}s:\tilde F^m\to F^m$.
Also, we set
\be
F_{a_1,\dots,a_n}(z_1,\dots,z_n,\lambda)
=\oplus_{m=0}^\infty F^m_{a_1,\dots,a_n}(z_1,\dots,z_n,\lambda),
\ee
and define an $\h$-module structure on
$F_{a_1,\dots,a_{n}}(z_1,\dots,z_{n},\lambda))$ by letting $h$ act by
\be
h|_{F^{m}_{a_1,\dots,a_{n}}(z_1,\dots,z_{n},\lambda)}=(\tsize\Sum_{i=1}^n
\Lambda_i-2m)\Id,\qquad a_i=\eta\Lambda_i.
\ee
\end{definition}

Clearly, $F^m_{a_{\sigma(1)},\dots,a_{\sigma(n)}}(z_{\sigma(1)},\dots,
z_{\sigma(n)})=F^{m}_{a_1,\dots,a_{n}}(z_1,\dots,z_{n},\lambda)$
for any permutation $\sigma\in S_n$.


It is shown in \cite{FTV} that
$F^{m}_{a_1,\dots,a_{n}}(z_1,\dots,z_{n},\lambda)$
is a finite-dimensional vector space of dimension
$\left(\begin{matrix}{n+m-1}\\{m}\end{matrix}\right)$.

\begin{example}
Let $n=1$. Then $F^m_a(z,\lambda)$ is a one-dimensional space spanned by
\begin{equation}\label{eb}
\omega_m(t_1,\dots,t_m,\lambda;z)=
\prod_{i<j}\frac{\theta(t_i-t_j)}
{\theta(t_i-t_j+2\eta)}\prod_{j=1}^{m}\frac
{\theta(\lambda+2\eta m+t_j-z-a)}
{\theta(t_j-z-a)}.
\end{equation}
\end{example}

\subsection{Tensor products}\label{sstp}
\begin{proposition}\label{p32} \cite{FTV}
Let $n=n'+n''$, $m=m'+m''$ be nonnegative integers and $a_1,\dots,a_n$,
$z_1,\dots,z_n$ be complex numbers. The formula
\be
k(t_1,\dots,t_{m})\,=\,
\frac1{m'!m''!}\;\Sym\biggl( f(t_{1},\dots,t_{m'})g(t_{m'+1},\dots,t_{m})
\prod_{
\begin{matrix} \scriptstyle{m'<j\leq m}\\
\scriptstyle{1\leq l\leq n'}\end{matrix}}
\frac{\theta(t_j-z_l+a_l)}{\theta(t_j-z_l-a_l)}
\biggr)\ee
correctly defines a linear map $\Phi:f\otimes g\mapsto k=\Phi(f\otimes g)$,
\ifTwelve
\bea
\oplus_{m'=0}^{m}
F^{m'}_{a_1,\dots,a_{n'}}(z_1,\dots,z_{n'},\lambda)\otimes
F^{m''}_{a_{n'+1},\dots,a_{n}}(z_{n'+1},\dots,z_{n},\lambda-2\nu)
\\
\to F^{m}_{a_1,\dots,a_{n}}(z_1,\dots,z_{n},\lambda),
\eea
\else
\bea
\oplus_{m'=0}^{m}F^{m'}_{a_1,\dots,a_{n'}}(z_1,\dots,z_{n'},\lambda)\otimes
F^{m''}_{a_{n'+1},\dots,a_{n}}(z_{n'+1},\dots,z_{n},\lambda-2\nu)\to
F^{m}_{a_1,\dots,a_{n}}(z_1,\dots,z_{n},\lambda),
\eea
\fi
where
$
\nu=a_{1}+\cdots +a_{n'}-2\eta m'.
$
For generic values of the parameters $z_j$, $\lambda$, the map $\Phi$ is
an isomorphism. Moreover, $\Phi$ is associative in the sense that, for any
three functions $f,g,h$, $\Phi(\Phi(f\otimes g)\otimes h)=
\Phi(f\otimes\Phi(g\otimes h))$, whenever defined.
\end{proposition}

By iterating this construction, we get for all $n\geq1$ a linear map $\Phi_n$,
defined recursively by $\Phi_1=\Id$, $\Phi_{n}=\Phi(\Phi_{n-1}\otimes\Id)$,
from
\be
\oplus_{m_1+\cdots+m_n=m}\otimes_{i=1}^n
F_{a_i}^{m_i}(z_i,\lambda-2\eta(\mu_{1}+\cdots+\mu_{i-1}))
\ee
to $F^m_{a_1,\dots,a_n}(z_1,\dots,z_n,\lambda)$, with $\mu_j=a_j/\eta-2m_j$,
$j=1,\dots,n$.

Let $V_\Lambda^*=\oplus_{j=0}^\infty\C e_j^*$ be the restricted dual of
the module $V_\Lambda=\oplus_{j=0}^\infty \C e_j$. It is spanned by the basis
$(e_j^*)$ dual to the basis $(e_j)$. We let $\h$ act on $V_\Lambda^*$ by
$he_j^*=(\Lambda-2j)e_j^*$. Then the map that sends $e^*_j$ to $\omega_j$
(see \Ref{eb}) defines an isomorphism of $\h$-modules
\be
\omega(z,\lambda):V_\Lambda^*\to F_a(z,\lambda), \qquad a=\eta\Lambda.
\ee
By composing this with the maps $\Phi$ of Proposition \ref{p32}, we obtain
homomorphisms (of $\h$-modules)
\be
\omega(z_1,\dots,z_n,\lambda):
V_{\Lambda_1}^*\otimes\cdots\otimes
V_{\Lambda_n}^*\to F_{a_1,\dots,a_n}(z_1,\dots,z_n,\lambda)
\ee
which are isomorphisms for generic values of $z_1,\dots,z_n,\lambda$.
The restriction of the map $\omega(z_1,\dots,z_n,\lambda)$ to
$\C e_{m_1}^*\otimes\cdots\otimes e_{m_n}^*$ is
\be
\Phi_n(\omega(z_1,\lambda) e^*_{m_1}
\otimes
\omega(z_{2},\lambda-2\eta\mu_1) e^*_{m_2}
\otimes
\cdots\otimes\omega(z_n,\lambda-2\eta(\mu_1+\cdots+\mu_{n-1}))
e^*_{m_n}) ,
\ee
where $\mu_j=\Lambda_j-2m_j$, $j=1,\dots,n$. For example, if $n=2$,
then $\omega(z_1,z_2,\lambda)$ sends $e^*_j\otimes e^*_k$ to
\be
\frac1{j!k!}\;\Sym\biggl(\omega_j(t_1,\dots,t_{j},\!\lambda;z_1)
\omega_k(t_{j+1},\dots,t_{j+k},\!\lambda-2a_1+4\eta j;z_2)\!
\prod_{i=j+1}^{j+k}\!\frac{\theta(t_i-z_1+a_1)}{\theta(t_i-z_1-a_1)}\biggr),
\ee
where $\{\omega_j(t_1,\dots,t_j,\lambda;z)\}$ is the basis \Ref{eb}
of $F_a(z,\lambda)$.

More generally, we have an explicit formula for the image of
$e_{m_1}^*\otimes\cdots\otimes e_{m_n}^*$, which we discuss next.

\subsection{A basis of $F_{a_1,\dots,a_n}(z_1,\dots,z_n)$}
The space $V_\Lambda$ comes with a basis $e_j$. Thus we have the natural basis
$e_{m_1}^*\otimes\cdots\otimes e_{m_n}^*$ of the tensor product of
$V^*_{\Lambda_i}$ in terms of the dual bases of the factors. The map
$\omega(z_1,\dots,z_n,\lambda)$ maps, for generic $z_i$, this basis to a basis
of $F_{a_1,\dots, a_n}(z_1,\dots,z_n,\lambda)$, which is an essential part of
our formulae for integral representations for solutions of the qKZB equations.

We give here an explicit formula for the basis vectors.
\begin{proposition}\label{pomegam} \cite{FTV}
Let $m\in\Z_{\geq 0}$, $\Lambda=(\Lambda_1,\dots,\Lambda_n)\in\C^n$, and let
$z=(z_1,\dots,z_n)\in\C^n$ be generic. Set $a_i=\eta\Lambda_i$. Let
\be
u(t_1,\dots,t_m)=\prod_{i<j}
\frac{\theta(t_i-t_j+2\eta)}
{\theta(t_i-t_j)}
\ee
Then, for generic $\lambda\in\C$, the functions
\be
\omega_{m_1,\dots,m_n}(t_1,\dots,t_m,\lambda;z)\,=\,\omega(z,\lambda)
\,e_{m_1}^*\otimes\cdots\otimes e_{m_n}^*
\ee
labeled by $\,{m_1,\dots,m_n\in\Z}\,$ with $\,{\sum_km_k=m}\,$ form a basis
of $F^m_a(z,\lambda)$ and are given by the explicit formula
\bean \label{ w.f }
\kern-20pt&& \omega_{m_1,\dots,m_n}(t_1,\dots,t_m,\lambda;z;\tau)\,=\,
u(t_1,\dots,t_m)^{-1}
\sum_{I_1,\dots,I_n}\prod_{l=1}^n\prod_{i\in I_l}\prod_{k=1}^{l-1}
\frac{\theta(t_i-z_k+a_k)}{\theta(t_i-z_k-a_k)}
\\
\kern-20pt&& {}\times\,\prod_{k<l}\prod_{i\in I_k,j\in I_l}
\frac{\theta(t_i-t_j+2\eta)}{\theta(t_i-t_j)}
\prod_{k=1}^{n}\prod_{j\in I_k}
\frac{\theta(\lambda\!+\!t_j\!-\!z_k\!-\!a_k\!+\!2\eta
m_k\!-\!2\eta\sum_{l=1}^{k-1}(\Lambda_l\!-\!2m_l))}{\theta(t_j-z_k-a_k)}\,.
\notag
\eean
The summation is over all $n$-tuples $I_1,\dots,I_n$ of disjoint subsets of
$\{1,\dots,m\}$ such that $I_j$ has $m_j$ elements, $1\leq j\leq n$.
\end{proposition}

We shall call the functions
$\omega_{m_1,\dots,m_n}(t_1,\dots,t_m,\lambda;z;\tau)$ the weight functions.

\subsection{$R$-matrices}
Let $a=\eta\Lambda$ and $b=\eta \Mu$ be complex numbers. Since
$F_{ab}(z,w,\lambda)$ coinsides with $F_{ba}(w,z,\lambda)$ by the symmetry
of the definition, we obtain a family of isomorphisms between
$V^*_\Lambda\otimes V^*_\Mu$ and $V^*_\Mu\otimes V^*_\Lambda $. The composition
of this family with the flip
$P:V^*_\Mu\otimes V^*_\Lambda\to V^*_\Lambda\otimes V^*_\Mu$,
$Pv\otimes w=w\otimes v$ gives a family of automorphisms of
$V^*_\Lambda\otimes V^*_\Mu$:

\vs\begin{definition}
Let $z,w,\lambda$ be such that
$\omega(z,w,\lambda):V^*_\Lambda \otimes V^*_\Mu\to F_{ab}(z,w,\lambda)$ is
invertible. The {\em $R$-matrix}
$R_{\Lambda,\Mu}(z,w,\lambda)\in\End_\h(V_\Lambda \otimes V_\Mu)$ is the dual
map to the composition $R^*_{\Lambda,\Mu}(z,w,\lambda)$:
\be
V^*_\Lambda\otimes V^*_\Mu
\buildrel{P}\over{\longrightarrow}
V^*_\Mu\otimes V^*_\Lambda
\buildrel{\omega(w,z,\lambda)}\over{\longrightarrow}
F_{ab}(z,w,\lambda)
\buildrel{\omega(z,w,\lambda)^{-1}}\over{\longrightarrow}
V^*_\Lambda \otimes V^*_\Mu,
\ee
where we identify canonically $V_\Lambda ^*\otimes V_\Mu^*$ with
$(V_\Lambda \otimes V_\Mu)^*$.
\end{definition}

Alternatively, the $R$-matrix $R_{\Lambda,\Mu}(z,w,\lambda)$ can be thought
of as the transition matrix expressing the basis
$\tilde\omega_{ij}=\omega(w,z,\lambda)e_j^*\otimes e_i^*$
of the space $F_{ab}(z,w,\lambda)$
in terms of the basis $\omega_{ij}=\omega(z,w,\lambda)e_i^*\otimes e_j^*$:
if $R_{\Lambda,\Mu}(z,w,\lambda)e_i\otimes e_j=
\sum_{kl}R_{ij}^{kl}e_k\otimes e_l$, then
\begin{equation} \label{R.c}
\tilde\omega_{kl}=\sum_{ij}R_{ij}^{kl}\omega_{ij}.
\end{equation}

\begin{lemma}\label{lpp} \cite{FTV}

\begin{enumerate}
\item[(i)] $R_{\Lambda, \Mu}(z,w,\lambda)$ is a meromorphic function
of $\Lambda,\Mu,z,w,\lambda$.
\item[(ii)] If $\Lambda$ is generic, then $R_{\Lambda,\Lambda}(z,w,\lambda)$
is regular at $z=w$ and $\lim_{z\to w}R_{\Lambda,\Lambda}(z,w,\lambda)=P$,
where $P$ is the flip $u\otimes v\mapsto v\otimes u$.
\item[(iii)]
$R_{\Lambda, \Mu}(z,w,\lambda)$ depends only on the difference $z-w$.
\end{enumerate}
\end{lemma}

Accordingly, we write $R_{\Lambda,\Mu}(z-w,\lambda)$ instead of
$R_{\Lambda,\Mu}(z,w,\lambda)$ in what follows.

The $R$-matrices satisfy the dynamical Yang-Baxter equation.

\begin{thm} \cite{FTV}
The matrices $R_{\Lambda,\Mu}(z,\lambda)$ obey I--III of Section \ref{sqkzb}.
\end{thm}

Let us now consider the case of positive integer weights. In this case
the $R$-matrices have invariant subspaces. If $\Lambda\in\Z_{\geq 0}$
we let $SV_\Lambda$ be the subspace of $V_\Lambda$ spanned by
$e_{\Lambda+1},\,e_{\Lambda+2},\dots$. The $\Lambda+1$-dimensional quotient
$V_\Lambda/SV_\Lambda$ will be denoted by $L_\Lambda$, and will be often
identified with $\oplus_{j=0}^m\C e_j$.

\begin{thm}\label{tfd} \cite{FTV}
Let $z,\eta,\lambda$ be generic and $\Lambda$, $\Mu\in\C$.
\begin{enumerate}
\item[(i)] If $\Lambda\in\Z_{\geq 0}$, then $R_{\Lambda,\Mu}(z,\lambda)$
preserves $SV_\Lambda\otimes V_\Mu$
\item[(ii)] If $\Mu\in\Z_{\geq 0}$, then $R_{\Lambda,\Mu}(z,\lambda)$
preserves $V_\Lambda\otimes SV_\Mu$
\item[(iii)] If $\Lambda\in\Z_{\geq 0}$ and $\Mu\in\Z_{\geq 0}$,
then $R_{\Lambda,\Mu}(z,\lambda)$ preserves
$SV_\Lambda\otimes V_\Mu+V_\Lambda\otimes SV_\Mu$.
\end{enumerate}
\end{thm}

In particular, if $\Lambda$ and/or $\Mu$ are nonnegative integers,
then $R_{\Lambda,\Mu}(z,\lambda)$ induces operators, still denoted by
$R_{\Lambda,\Mu}(z,\lambda)$, on the quotients $L_\Lambda\otimes V_\Mu$,
$V_\Lambda\otimes L_\Mu$ and/or $L_\Lambda\otimes L_\Mu$.
They obey the dynamical Yang--Baxter equation.

\subsection{Example}\label{ssexa}
We give an example of computation of matrix elements of the $R$-matrix
$R_{\Lambda,\Mu}(z-w,\lambda)$, assuming that the parameters are generic.

The $R$-matrix is calculated as the transition matrix relating two bases of
$F_{ab}(z,w,\lambda)$: let
\be
\tilde\omega_{ij}=\omega(w,z,\lambda)\,e_j^*\otimes e_i^*,\qquad
\omega_{ij}=\omega(z,w,\lambda)\,e_i^*\otimes e_j^*.
\ee
The matrix elements of $R$ with respect to the basis $e_j\otimes e_k$ are
given by $\tilde\omega_{kl}=\sum_{ij}R_{ij}^{kl}\omega_{ij}$.
The $R$-matrix preserves the weight spaces
\be
(V_\Lambda\otimes V_{\Mu})[\Lambda+\Mu-2m]=\oplus_{j=0}^m\C e_j\otimes
e_{m-j} ,
\ee
and we may consider the problem of computing the matrix elements of
the $R$-matrix separately on each weight space. Without loss of generality
we assume that $w=0$.

Let $m=0$. Then the weight space is spanned by $e_0\otimes e_0$ and
$\omega_{00}=\tilde\omega_{00}=1$. Therefore $R_{00}^{00}=1$.

Let $m=1$. Then the basis elements are functions of one variable $t=t_1$
and we have (with $a=\eta\Lambda$, $b=\eta\Mu$)
\be
\omega_{01}(t)=\frac{\theta(\lambda+2\eta+t-2a-b)\theta(t-z+a)}
{\theta(t-b)\theta(t-z-a)}\;,
\qquad
\omega_{10}(t)=\frac{\theta(\lambda+2\eta+t-z-a)}{\theta(t-z-a)}\;,
\ee
and
\be
\tilde\omega_{01}(t)=\frac{\theta(\lambda+2\eta+t-b)}{\theta(t-b)}\;,
\qquad
\tilde\omega_{10}(t)=\frac{\theta(\lambda+2\eta+t-z-2b-a)\theta(t+b)}
{\theta(t-z-a)\theta(t-b)}\;.
\ee

\begin{proposition}\label{Rijkl} \cite{FTV}

Let the matrix elements of the $R$-matrix $R_{\Lambda,\Mu}(z,\lambda)$ be
defined by
\be
R_{\Lambda,\Mu}(z,\lambda)e_i\otimes e_j=\tsize\sum_{kl}
R_{ij}^{kl}
e_k\otimes e_l.
\ee
Then
\bea
R^{00}_{00}\,=&&\kern-14pt 1\,,
\\
R^{01}_{01}\,=&&\kern-14pt
\frac{\theta(z+\eta\Lambda-\eta\Mu)\theta(\lambda+2\eta)}
{\theta(z-\eta\Lambda-\eta\Mu)\theta(\lambda+2\eta(1-\Lambda))}\;,
\\
R^{01}_{10}\,=&&\kern-14pt
-\frac{\theta(\lambda+2\eta+z-\eta\Lambda-\eta\Mu)\theta(2\eta\Lambda)}
{\theta(z-\eta\Lambda-\eta\Mu)\theta(\lambda+2\eta(1-\Lambda))}\;,
\\
R^{10}_{01}\,=&&\kern-14pt
-\frac{\theta(\lambda+2\eta-z-\eta\Lambda-\eta\Mu)\theta(2\eta\Mu)}
{\theta(z-\eta\Lambda-\eta\Mu)\theta(\lambda+2\eta(1-\Lambda))}\;,
\\
R^{10}_{10}\,=&&\kern-14pt
\frac{\theta(z+\eta\Mu-\eta\Lambda)\theta(\lambda+2\eta(1-\Lambda-\Mu))}
{\theta(z-\eta\Lambda-\eta\Mu)\theta(\lambda+2\eta(1-\Lambda))}\;.
\eea
\end{proposition}

\subsection{Evaluation Verma modules and their tensor products}
\label{ssevm}
Here we explain the relation between the geometric construction of
tensor products and $R$-matrices and the representation theory of
$E_{\tau,\eta}(sl_2)$ \cite{FV1}.

Recall the definition of a representation of $E_{\tau,\eta}(sl_2)$:
let $\h$ act on $\C^2$ via $h={\mathrm {diag}}(1,-1)$. A representation of
$E_{\tau,\eta}(sl_2)$ is an $\h$-module $W$ with diagonalizable action of $h$
and finite-dimensional eigenspaces, together with an operator
$L(z,\lambda)\in\End(\C^2\otimes W)$ (the ``$L$-operator''),
commuting with $h^{(1)}+h^{(2)}$, and obeying the relations
\ifTwelve
\bea
R^{(12)}(z-w,\lambda-2\eta h^{(3)}) &&\kern-20pt\, L^{(13)}(z,\lambda)\,
L^{(23)}(w,\lambda-2\eta h^{(1)})
\\[2pt]
{}={}&&\kern-20pt\, L^{(23)}(w,\lambda)\,
L^{(13)}(z,\lambda-2\eta h^{(2)})\,R^{(12)}(z-w,\lambda)
\eea
\else
\bea
R^{(12)}(z-w,\lambda-2\eta h^{(3)})\,L^{(13)}(z,\lambda)\,
L^{(23)}(w,\lambda-2\eta h^{(1)})\,=\,L^{(23)}(w,\lambda)\,
L^{(13)}(z,\lambda-2\eta h^{(2)})\,R^{(12)}(z-w,\lambda)
\eea
\fi
in $\End(\C^2\otimes\C^2\otimes W)$. The {\em fundamental $R$-matrix}
$R(z,\lambda)\in\End(\C^2\otimes\C^2)$ is the following solution of the
dynamical Yang--Baxter equation: let $e_0$, $e_1$ be the standard basis of
$\C^2$, then with respect to the basis $e_0\otimes e_0$, $e_0\otimes e_1$,
$e_1\otimes e_0$, $e_1\otimes e_1$ of $\C^2\otimes\C^2$,
\be
R(z,\lambda)=
\left(\begin{matrix}
1 & 0 & 0 & 0\\
0 & \alpha(z,\lambda) & \beta(z,\lambda) & 0\\
0 & \beta(z,-\lambda) & \alpha(z,-\lambda) & 0\\
0 & 0 & 0 & 1
\end{matrix}\right).
\ee
where
\be
\alpha(z,\lambda)=
\frac{\theta(\lambda+2\eta)\theta(z)}{\theta(\lambda)\theta(z-2\eta)},\qquad
\beta(z,\lambda)=
-\frac{\theta(\lambda+z)\theta(2\eta)}{\theta(\lambda)\theta(z-2\eta)}.
\ee
To discuss representation theory, it is convenient to think of
$L(z,\lambda)\in\End(\C^2\otimes W)$ as a two by two matrix with entries
$a(z,\lambda)$, $b(z,\lambda)$, $c(z,\lambda)$, $d(z,\lambda)$ in $\End(W)$.
In \cite{FV1} we wrote explicitly the relations that these four operators must
satisfy, and defined a class of representations, the evaluation Verma modules,
by giving explicitly the action of these four operators on basis vectors.
These formulae can be obtained from the geometric construction.

\begin{thm}
\label{tfrm} \cite{FTV}
Let us identify the two-dimensional space $L_1=V_{\Lambda=1}/SV_{\Lambda=1}$
with $\C^2$ via the basis $e_0,e_1$. Then the $R$-matrix
$R_{1,1}(z,\lambda)\in\End(L_1\otimes L_1)$ coincides with
the fundamental $R$-matrix.
\end{thm}

\begin{corollary}
For any $w,\Mu\in\C$, the $\h$-module $V_\Mu$ together with the operator
$L(z,\lambda)=R_{1,\Mu}(z-w,\lambda)\in\End(L_1\otimes V_\Mu)$ defines
a representation of $E_{\tau,\eta}(sl_2)$.
\end{corollary}

This representation is called in \cite{FV1} the evaluation Verma module with
evaluation point $w$ and highest weight $\Mu$. It is denoted by $V_\Mu(w)$.
The matrix elements of $L(z,\lambda)$ are given explicitly in Theorem 3 of
\cite{FV1} in terms of the action of $a(z,\lambda)$,\dots, $d(z,\lambda)$.
In the notation of Proposition \ref{Rijkl}, this result amounts to
the following formulae for the matrix elements $R_{ij}^{kl}$
of $R_{1,\Lambda}(z,\lambda)\in \End(L_1\otimes V_\Lambda)$.
\bea
R_{0k}^{0k}\,=&&\kern-14pt
\frac{\theta(z-(\Lambda+1-2k)\eta)}{\theta(z-(\Lambda+1)\eta)}
\frac{\theta(\lambda+2k\eta)}{\theta(\lambda)},
\\
R_{1k}^{0,k+1}\,=&&\kern-14pt
-\,\frac{\theta(\lambda+z-(\Lambda-1-2k)\eta)}{\theta(z-(\Lambda+1)\eta)}
\frac{\theta(2\eta)}{\theta(\lambda)},
\\
R_{0k}^{1,k-1}\,=&&\kern-14pt
-\,\frac{\theta(\lambda-z-(\Lambda+1-2k)\eta)}{\theta(z-(\Lambda+1)\eta)}
\frac{\theta(2(\Lambda+1-k)\eta)}{\theta(\lambda)}
\frac{\theta(2k\eta)}{\theta(2\eta)},
\\
R_{1k}^{1k}\,=&&\kern-14pt
\frac{\theta(z-(-\Lambda+1+2k)\eta)}{\theta(z-(\Lambda+1)\eta)}
\frac{\theta(\lambda-2(\Lambda-k)\eta)}{\theta(\lambda)}.
\eea
Moreover, the tensor product construction of \ref{sstp} is related to
the tensor product of representations of the elliptic quantum group.
Recall that if $W_1$, $W_2$ are representations of the elliptic quantum group
with $L$-operators $L_1(z,\lambda)$, $L_2(z,\lambda)$, then
their tensor product $W=W_1\otimes W_2$ with $L$-operator
\be
L(z,\lambda)=L_1(z,\lambda-2\eta h^{(3)})^{(12)}
L_2(z,\lambda)^{(13)}\in\End(\C^2\otimes W)
\ee
is also a representation of the elliptic quantum group.

\begin{thm}\cite{FTV}

Let $\Lambda_1,\dots,\Lambda_n\in\C$ and $z_1,\dots,z_n$ be generic complex
numbers. Let $V=V_{\Lambda_1}\otimes\cdots\otimes V_{\Lambda_n}$ and
$L(z,\lambda)\in\End(V_{\Lambda=1}\otimes V)$ be defined by the relation
\be
\omega(z,z_1,\dots,z_n,\lambda)L(z,\lambda)^*=
\omega(z_1,\dots,z_n,z,\lambda)P
\ee
in $\End((V_1\otimes V)^*)=\End(V_1^*\otimes V^*)$, where
$Pv_1\otimes v=v\otimes v_1$, if $v_1\in V_1^*$, $v\in V^*$.
Then $L(z,\lambda)$ is well-defined as an endomorphism of the quotient
$L_1\otimes V=\C^2\otimes V$, and defines a structure of a representation of
$E_{\tau,\eta}(sl_2)$ on $V$. This representation is isomorphic to the tensor
product of evaluation Verma modules
\be
V_{\Lambda_n}(z_n)\otimes\cdots\otimes V_{\Lambda_1}(z_1),
\ee
with the isomorphism
$u_1\otimes\cdots\otimes u_n\mapsto u_n\otimes\cdots\otimes u_1$.
\end{thm}

Finally, the dynamical Yang--Baxter equation in
$L_1\otimes V_\Lambda\otimes V_\Mu$ can be stated as saying
that $R_{\Lambda,\Mu}(z-w,\lambda)P$ is an isomorphism from
$V_\Mu(w)\otimes V_\Lambda(z)$ to $V_\Lambda(z)\otimes V_\Mu(w)$,
see \cite{FV1}. Therefore, by uniqueness, the $R$-matrices constructed
in this section coincide with the solutions of the dynamical Yang--Baxter
equation described in Section 13 of \cite{FV1}.

\section{Transformation of qKZB with respect to shifts
of $z_k$ by $\tau$ and $1$}\label{transf}

\subsection{ Transformation of weight functions}

Let $\omega_{m_1,\dots,m_n}(t,\lambda;z_1,\ldots,z_n;\tau)$ be the weight
functions defined by \Ref{ w.f }. The next Proposition describes the
transformation properties of the weight functions with respect to shifts
by $\tau$ of variables $z_1,\ldots,z_n$.

\begin{proposition}\label{trans}
For any $k$, we have
\be
\omega_{m_1,\dots,m_n}(t,\lambda;z_1,\ldots,z_k+\tau,\ldots,z_n;\tau)
=e^{b_k}\omega_{m_1,\dots,m_n}(t,\lambda;z_1,\ldots,z_k,\ldots,z_n;\tau)
\ee
where
\be
b_k=2\pi im_k (\lambda + 2\eta m_k - 2\sum_{l=1}^{k-1}
(a_l - 2\eta m_l)) + 4\pi i a_k \sum_{l=k+1}^n m_{l} .
\ee
\end{proposition}

The Proposition follows from transformation properties of $\theta(t)$.

We reformulate the Proposition. Introduce a function
\begin{equation}\label{alpha}
\al(\la)\,=\, \exp (-\pi i \la^2/4\eta)\,.
\end{equation}
Use the notation $h^{(j)} = a_j/\eta -2m_j$. Set
\begin{equation}\label{DD}
D_{k}(\la)\,=\,
\frac{\al(\la\,-\,2\eta\,\sum_{l=1}^{k}h^{(l)})}{
\al(\la\,-\,2\eta\,\sum_{l=1}^{k-1}h^{(l)})}\,
e^{\pi i\eta\Lambda_k(\sum_{l=1}^{k-1}\Lambda_l-\sum_{l=k+1}^{n}\Lambda_l)}.
\end{equation}
Then
\begin{equation}\label{b}
e^{b_k}\,=\,D_k^{-1}\,e^{\pi i\Lambda_k(\la+2\eta m) }\,
e^{-\pi i\eta\Lambda_k\sum_{l=1}^nh^{(l)}}\,,
\end{equation}
where $m=m_1+\ldots+m_n$.

\subsection{ Transformations of $R$-matrices}

Consider the tensor product of evaluation Verma modules over the elliptic
quantum group $E_{\tau, \eta}(sl_2)$, $V_\Lambda(z)\otimes V_\Mu(w)$, and its
$R$-matrix $R_{\Lambda,\Mu}(z-w,\la)\in\End (V_\Lambda\otimes V_\Mu)$ defined
in Section \ref{ssevm}. Let $h^{(j)}$ denote the operator $h$ acting in
the $j$-th factor of $V_\Lambda(z)\otimes V_\Mu(w)$.

\begin{proposition}\label{transR}
\bean\label{t.r}
R_{\Lambda,\Mu}(z+\tau,\la)\,=\,e^{-2\pi i\Lambda\Mu}\,\,
\frac{\al(\la\,-\,2\eta h^{(2)})}{\al(\la\,-\,2\eta\,(h^{(1)}+ h^{(2)}))}\,
R_{\Lambda,\Mu}(z,\la)\,\,\frac{\al(\la\,-\,2\eta h^{(1)})}{\al(\la)}\,.
\eean

\end{proposition}

\begin{proof} Use formula \Ref{R.c},
\be
\tilde\omega_{kl}(z,w)\,=\,\sum_{ij}\,R_{ij}^{kl}(z-w,\la)\,\omega_{ij}(z,w).
\ee
Then
\be
\tilde\omega_{kl}(z+\tau,w)\,=\,
\sum_{ij}\,R_{ij}^{kl}(z-w+\tau,\la)\,\omega_{ij}(z+\tau,w).
\ee
We have $\tilde\omega_{kl}(z+\tau,w)=e^{\tilde b_1(kl)}\tilde\omega_{kl}(z,w)$
and $\omega_{ij}(z+\tau,w)=e^{b_1(ij)}\omega_{ij}(z,w)$ where $e^{\tilde b_1}$
and $e^{b_1}$ are given by \Ref{b}. Hence
\be
R_{ij}^{kl}(z-w+\tau,\la)\,=\,
e^{\tilde b_1(kl)}\,R_{ij}^{kl}(z-w,\la)\,e^{-b_1(ij)}.
\ee
This proves the Proposition.
\end{proof}

\begin{proposition}\label{l-trans}
\bean\label{L.trans}
R_{\Lambda,\Mu}(z-w,\la+\tau)\,=\,
e^{\pi i(h^{(1)}(-z-\eta\Mu)+h^{(2)}(-w+\eta\Lambda))}\,
R_{\Lambda,\Mu}(z-w,\la)\,
e^{\pi i(h^{(1)}(z-\eta\Mu)+h^{(2)}(w+\eta\Lambda))}\,.
\notag
\eean
\end{proposition}
The proof is analogous to the proof of Proposition \ref{transR}.
$\square$

\begin{proposition}\label{tau-trans}
\bean\label{tau.trans}
R_{\Lambda,\Mu}(z+1,\la, \tau)\,=\,R_{\Lambda,\Mu}(z,\la+1, \tau)\,=\,
R_{\Lambda,\Mu}(z,\la, \tau +1)\,=\,R_{\Lambda,\Mu}(z,\la, \tau)\,.
\notag
\eean
\end{proposition}
The proof easily follows from the formulae $\theta(t+1,\tau)=-\theta(t,\tau)$
and $\theta(t,\tau + 1)=e^{\pi i/4}\theta(t,\tau)$.

\subsection{ Transformation of qKZB equations with respect to shifts
$z_k\to z_k+\tau$}

Consider the qKZB equations defined by \Ref{qKZB},
\begin{equation}
\Psi(z_1,\dots,z_j+p,\dots,z_n,\la)\,=\,K_j(z_1,\dots,z_n,\la)\,
\Psi(z_1,\dots,z_n,\la)\,,\qquad j=1,\dots,n,
\notag
\end{equation}
where $\Psi(z_1,\dots,z_n,\la)$ is a function with values in
$V=V_{\Lambda_1}\otimes\dots\otimes V_{\Lambda_n}$ and
$K_j(z_1,\dots,z_n,\la)\in\End(V)$.

For any $k=1,\dots,n$, introduce a linear operator
$B_k(z_1,\dots,z_n,\la)\in\End(V[0])$ by
\ifTwelve
\bean \label{Bk}
\kern-18pt&& B_k(z_1,\dots,z_n,\la)\,={}
\\
\kern-18pt&& e^{2\pi i\eta\sum_{l,\,l\neq k}(z_l-z_k)\Lambda_l\Lambda_k/p}
\,\,\frac{\al(\la\,-\,2\eta\,\sum_{l=1}^{k}h^{(l)})}
{\al(\la\,-\,2\eta\,\sum_{l=1}^{k-1}h^{(l)})}\,
e^{\pi i\eta\Lambda_k(\sum_{l=1}^{k-1}\Lambda_l-\sum_{l=k+1}^{n}\Lambda_l)}.
\notag
\eean
\else
\bean \label{Bk}
B_k(z_1,\dots,z_n,\la)\,=
\,e^{2\pi i\eta\sum_{l,\,l\neq k}(z_l-z_k)\Lambda_l\Lambda_k/p}
\,\,\frac{\al(\la\,-\,2\eta\,\sum_{l=1}^{k}h^{(l)})}
{\al(\la\,-\,2\eta\,\sum_{l=1}^{k-1}h^{(l)})}\,
e^{\pi i\eta\Lambda_k(\sum_{l=1}^{k-1}\Lambda_l-\sum_{l=k+1}^{n}\Lambda_l)}.
\eean
\fi

\begin{thm}
\label{tr.qkzb}
For any $j$ and $k$, we have
\bean\label{tr.op} 
\kern-18pt&& K_j(z_1,\dots,z_k+\tau,\dots,z_n,\la)\,={}
\\[2pt]
\kern-18pt&& B_k(z_1,\dots,z_j+p,\dots,z_n,\la)^{-1}\,K_j(z_1,\dots,z_n,\la)\,
B_k(z_1,\dots,z_j,\dots,z_n,\la)\,.
\notag
\eean

\end{thm}
Notice that the last exponential in \Ref{Bk} is not essential for
the Theorem and is introduced for later purposes.

\begin{proof}
There are three cases: $k>j,\, k=j, \,k<j.$ We prove the Theorem for $k>j$.
The other two cases are proved similarly.
\bean
K_j(z_1,\dots,z_k+\tau,\dots,z_n,\la)\,=\,R_{j,j-1}(z_j\!-\!z_{j-1}+p)\cdots
R_{j,1}(z_j\!-\!z_1+p)\Gamma_j R_{j,n}(z_j\!-\!z_n)\cdots
\notag
\\
R_{j,k}(z_j\!-\!z_k\,-\tau)\cdots R_{j,j+1}(z_j\!-\!z_{j+1})\,=\,
R_{j,j-1}(z_j\!-\!z_{j-1}+p)\cdots R_{j,1}(z_j\!-\!z_1+p)
\Gamma_j R_{j,n}(z_j\!-\!z_n)\cdots
\notag
\\
e^{2\pi i \eta \La_j\La_k}
\frac{\al(\la\,-\,2\eta\,\sum_{l=1}^k h^{(l)})}
{\al(\la\,-\,2\eta\,\sum_{l=1,\,l\neq j}^k h^{(l)})}
R_{j,k}(z_j\!-\!z_k)
\frac{\al(\la\,-\,2\eta\,\sum_{l=1,\,l\neq j}^{k-1}h^{(l)})}
{\al(\la\,-\,2\eta\,\sum_{l=1}^{k-1}h^{(l)})}
\cdots R_{j,j+1}(z_j\!-\!z_{j+1})\,=\,
\notag
\\
e^{2\pi i \eta \La_j\La_k}R_{j,j-1}(z_j\!-\!z_{j-1}+p)\cdots
R_{j,1}(z_j\!-\!z_1+p)\Gamma_j R_{j,n}(z_j\!-\!z_n)\cdots
\frac{\al(\la\,-\,2\eta\,\sum_{l=1,\,l\neq j}^{k-1}h^{(l)})}
{\al(\la\,-\,2\eta\,\sum_{l=1,\,l\neq j}^k h^{(l)})}
\notag
\\
R_{j,k}(z_j\!-\!z_k)
\frac{\al(\la\,-\,2\eta\,\sum_{l=1}^k h^{(l)})}
{\al(\la\,-\,2\eta\,\sum_{l=1}^{k-1}h^{(l)})}
\cdots R_{j,j+1}(z_j\!-\!z_{j+1})\,=\,e^{2\pi i \eta \La_j\La_k}
R_{j,j-1}(z_j\!-\!z_{j-1}+p)\cdots
\notag
\\
R_{j,1}(z_j\!-\!z_1+p)\Gamma_j
\frac{\al(\la\,-\,2\eta\,\sum_{l=1,\,l\neq j}^{k-1}h^{(l)})}
{\al(\la\,-\,2\eta\,\sum_{l=1,\,l\neq j}^k h^{(l)})}
R_{j,n}(z_j\!-\!z_n)\cdots R_{j,k}(z_j\!-\!z_k)\cdots
R_{j,j+1}(z_j\!-\!z_{j+1})
\notag
\\
\frac{\al(\la\,-\,2\eta\,\sum_{l=1}^k h^{(l)})}
{\al(\la\,-\,2\eta\,\sum_{l=1}^{k-1}h^{(l)})}\,=\,e^{2\pi i \eta \La_j\La_k}
\frac{\al(\la\,-\,2\eta\,\sum_{l=1}^{k-1}h^{(l)})}
{\al(\la\,-\,2\eta\,\sum_{l=1}^k h^{(l)})}
R_{j,j-1}(z_j\!-\!z_{j-1}+p)\cdots
\notag
\\
R_{j,1}(z_j\!-\!z_1+p)\Gamma_j R_{j,n}(z_j\!-\!z_n)\cdots
R_{j,k}(z_j\!-\!z_k)\cdots R_{j,j+1}(z_j\!-\!z_{j+1})
\frac{\al(\la\,-\,2\eta\,\sum_{l=1}^k h^{(l)})}
{\al(\la\,-\,2\eta\,\sum_{l=1}^{k-1}h^{(l)})}\,=\,
\notag
\\
e^{2\pi i \eta \La_j\La_k}
\frac{\al(\la\,-\,2\eta\,\sum_{l=1}^{k-1}h^{(l)})}
{\al(\la\,-\,2\eta\,\sum_{l=1}^k h^{(l)})}
K_j(z_1,\dots,z_n,\la)
\frac{\al(\la\,-\,2\eta\,\sum_{l=1}^k h^{(l)})}
{\al(\la\,-\,2\eta\,\sum_{l=1}^{k-1}h^{(l)})}.
\notag
\eean
This proves the Theorem for $k>j$.
\end{proof}

\begin{thm}
\label{mon.shi}
Let $\Psi(z_1,\dots,z_n,\la)$ be a solution of the qKZB equations \Ref{qKZB}.
Then for any $k$, the function
$B_k(z_1,\dots,z_n,\la)\Psi(z_1,\dots,z_k+\tau,\dots,z_n,\la)$
is a new solution of the same equations.
\end{thm}

\begin{proof}
For any $j$, we have
\bean
B_k(z_1,\dots,z_j+p,\dots,z_n,\la)\,
\Psi(z_1,\dots,z_k+\tau,\dots,z_j+p,\dots,z_n,\la)\,=\,
\notag
\\
B_k(z_1,\dots,z_j+p,\dots,z_n,\la)\,K_j(z_1,\dots,z_k+\tau,\dots,z_n,\la)\,
\Psi(z_1,\dots,z_k+\tau,\dots,z_j,\dots,z_n,\la)\,=\,
\notag
\\
B_k(z_1,\dots,z_j+p,\dots,z_n,\la)\,B_k(z_1,\dots,z_j+p,\dots,z_n,\la)^{-1}\,
K_j(z_1,\dots,z_n,\la)\,B_k(z_1,\dots,z_n,\la)\,
\notag
\\
\Psi(z_1,\dots,z_k+\tau,\dots,z_n,\la)\,=\,
K_j(z_1,\dots,z_n,\la)\,B_k(z_1,\dots,z_n,\la)\,
\Psi(z_1,\dots,z_k+\tau,z_n,\la)\,.
\notag
\eean

\end{proof}

\subsection{ Transformation of qKZB equations with respect to shifts
$z_k \to z_k+1$}

\begin{proposition}
\label{sh.1}
For any $j$, we have
$$
K_j(z_1,\dots,z_k+1,\dots,z_n, \,\la)\,=\,
K_j(z_1,\dots,z_k ,\dots,z_n, \,\la)\,.
$$
\end{proposition}
The Proposition follows from Proposition \ref{tau-trans}.

\begin{corollary}\label{shifts-1}
\label{mon.shi.1}
Let $\Psi(z_1,\dots,z_n,\la)$ be a solution of the qKZB equations \Ref{qKZB}.
Then for any $k$, the function $\Psi(z_1,\dots,z_k+1,\dots,z_n,\la)$ is a new
solution of the same equations.
\end{corollary}

\subsection{ The qKZB equations as equations on a torus }

Propositions \ref{mon.shi} and \ref{mon.shi.1} mean that we can define a vector
bundle over the torus $\C^n\,/\,\Z^n+\tau\Z^n$ whose fiber is the space of
functions of $\la$ with values in $V[0]$. The identification of points of
the base, $(z_1,\dots,z_n)\to(z_1,\dots,z_k+\tau,\dots,z_n)$, corresponds
to the identification of the fibers defined by $v\to B_k(z_1,\dots,z_n)v$.
The identification of points of the base,
$(z_1,\dots,z_n)\to (z_1,\dots,z_k+1,\dots,z_n)$, corresponds to
the identification of the fibers defined by $v\to v$. The identifications
commute with the qKZB equations. Hence the qKZB equations induce a system of
difference equations on the torus with values in this bundle (a flat discrete
connection). A solution of the initial qKZB equations defines a multivalued
solution of the difference equations on the torus. Now we can ask a question
about the monodromy of solutions of the equations on the torus.
We shall address this problem in Section \ref{Monodromy}

\section{Formal solutions of the qKZB equations}\label{ssqkzb}
In this section we fix $\tau,\eta,p,\Lambda_1,\dots,\Lambda_n$, and set
$a_i=\eta \Lambda_i$.

\subsection{Formal integral solutions}

By a formal Jackson integral solution of the qKZB equations we mean
an expression
\be
\Psi(z_1,\dots,z_n,\lambda)=
\int f(z_1,\dots,z_n,t_1,\dots,t_m,\lambda)Dt_1\cdots Dt_m,
\ee
where $f$ takes its values in $V[0]$, which obeys the qKZB equations \Ref{qKZB}
if we formally use the rule that the ``integral'' $\int$ is invariant under
translations of the variables $t_i$ by $p$. In other words,
$f(z_1,\dots,z_n,t_1,\dots,t_m ,\lambda)$ obeys the qKZB equations in
the variables $z_i$ up to terms of the form
$g(\dots,t_i+p,\dots)-g(\dots,t_i,\dots)$.

\medskip
\begin{definition}
A function $\Phi_{a}(t)$ depending on a complex parameter $a$, such that
\be\label{ph}
\Phi_a(t+p)=\frac{\theta(t+a)}{\theta(t-a)}\Phi_a(t)
\ee
is called a (one-variable) phase function.
\end{definition}

We assume that $p$ has positive imaginary part, and set $r=e^{2\pi ip}$,
$q=e^{2\pi i\tau}$. Then the convergent infinite product
\begin{equation}\label{phase1}
\Omega_a(t):=\Omega_a(t,\tau,p)=\prod_{j=0}^{\infty}\prod_{k=0}^{\infty}
\frac{(1-r^jq^ke^{2\pi i(t-a)})(1-r^{j+1}q^{k+1}e^{-2\pi i(t+a)})}
{(1-r^jq^ke^{2\pi i(t+a)})(1-r^{j+1}q^{k+1}e^{-2\pi i(t-a)})}\,,
\end{equation}
defines a phase function
\begin{equation}
\Phi_a(t)=e^{-2\pi iat/p}\Omega_a(t)
\notag
\end{equation}
and any other phase function is obtained from this one by multiplication
by a $p$-periodic function.

Given a one-variable phase function $\Phi_a(t)$, we define with our data
an $m$-variable phase function
\begin{equation}\label{phasem}
\Phi(t_1,\dots,t_m,z_1,\dots,z_n)=\prod_{j=1}^m\prod_{l=1}^n\Phi_{a_l}(t_j-z_l)
\prod_{1\leq i<j\leq m}\Phi_{-2\eta}(t_i-t_j).
\end{equation}

\begin{thm}\label{irs0} \cite {FTV}.
Let $\Phi_a(t)$ be a phase function, and let $\Phi$ be the corresponding
$m$-variable phase function \Ref{phasem}. For any entire function $\xi$
of one variable, let
\be\label{totalph}
\psi^\xi(t,z,\lambda)=\xi(p\lambda\!-\!
{\tsize\Sum_{l=1}^n2a_lz_l\!+\!4\eta\Sum_{j=1}^mt_j}))
\sum_{j_1+\cdots+j_n=m}\omega_{j_1,\dots,j_n}(t_1,\dots,t_m,\lambda)\,
e_{j_1}\!\otimes\!\cdots\!\otimes\! e_{j_n}.
\ee
Then
\be
\Psi(z_1,\dots,z_n,\lambda)=\int\Phi(t_1,\dots,t_m,z_1,\dots,z_n)
\psi^\xi(t_1,\dots,t_m,z_1,\dots,z_n,\lambda)Dt_1\cdots Dt_m
\ee
is a formal Jackson integral solution of the qKZB equations.
\end{thm}

The proof of the theorem is based on the transformation properties
of the $m$-variable phase function $\Phi(t,z)$ with respect to shifts
$t_j\mapsto t_j+p,\,z_l\mapsto z_l+p$ and does not use the explicit form
of the function $\Phi(t,z)$, see \cite {FTV}.

To obtain solutions from formal Jackson integral solutions, we need to find
{\em cycles}, linear forms on the space of functions of $t_1,\dots,t_m$ that
are invariant under translations $t_i\mapsto t_i+p$. To this end we need a
stronger version of the preceding theorem, Theorem \ref{irs} below gives us
a space of functions on which our cycles should be defined.

Let $\Phi$ be the phase function \Ref{phasem} and let $a=(a_1,\dots,a_n)$,
$z=(z_1,\dots,z_n)$. We assume, as usual, that $\sum a_i=2\eta m$,
$m\in\Z_{\geq 0}$. For any entire function $\xi$ of one variable, let
$E^0_a(z;\xi)$, be the space spanned by the functions of $t\in\C^m$ of the form
\be
\tsize \Phi(t,z)\xi(p\lambda-\Sum_{k=1}^n2a_kz_k+4\eta\Sum_{j=1}^mt_j)
f(t,z),
\ee
where $f(t,z)$, viewed as a function of $t=(t_1,\dots,t_m)$ belongs to
$\tilde F^m_a(z,\lambda)$ (see \ref{sssym}) for some $\lambda$.
All components of our integrand belong to this space.

Let $E_a(z;\xi)$, the space of cocycles, be the space spanned by functions
of the form $g(t+p\alpha)$, where $g\in E^0_a(z;\xi)$ and $\alpha\in\Z^m$.
By construction, $E_a(z;\xi)$ is invariant under translations of the arguments
$t_i$ by $p$. We define the space of coboundaries $DE_a(z;\xi)$ to be
the subspace of $E_a(z;\xi)$ spanned by functions of the form
$f(\dots,t_j+p,\dots)-f(\dots,t_j,\dots)$, $f\in E_a(z;\xi)$.

\begin{proposition}\label{pE(z)} \cite {FTV}.
$E_a(z_1,\dots,z_n)=E_a(z_1,\dots,z_j+p,\dots,z_n)$ for $j=1,\dots,n$
\end{proposition}

Moreover, we have the following result.

\begin{thm}\label{irs} \cite {FTV}.
Let us write the qKZB equations as $\Psi(\dots,z_j+p,\dots)=K_j(z)\Psi(z)$.
Then, for any entire function $\xi$, the integrand $\Phi(t,z)\psi^\xi(t,z)$
of Theorem \ref{irs0}, viewed as a function of $t\in\C^m$, belongs to
$E_a(z;\xi)\otimes V[0]$ for all $z\in\C^n$. It obeys the equations
\be
\Psi(t,\dots,z_j+p,\dots)=K_j(z)\Psi(t,z)\mod DE_a(z;\xi)\otimes V[0],
\qquad j=1,\dots,n.
\ee
In other words, $\Psi(z,t)$ solves the qKZB equations
in the cohomology $(E_a(z;\xi)/DE_a(z;\xi))\otimes V[0]$.
\end{thm}

To obtain solutions from these formal solutions, one should find horizontal
families of cycles, i.e., linear functions $\gamma(z)$ on $E_a(z;\xi)$
vanishing on $DE_a(z;\xi)$, and such that $\gamma(z+p\alpha)=\gamma(z)$
for all $\alpha\in\Z^n$. This problem will be addressed in the next section.

\subsection{New formal solutions}

Let $a=(a_1,\ldots,a_n),\,z=(z_1,\ldots,z_n),\,t=(t_1,\ldots,t_m)$. Assume
that $\sum a_i=2\eta m$ and the step $p$ of the qKZB equations has a positive
imaginary part. In this section we shall construct a finite set of new formal
integral solutions to the qKZB equations \Ref{qKZB}. The formal solutions are
labelled by an index $M=(m_1,\ldots,m_n)$ where non-negative integers $m_i$
satisfy $m_1+\ldots+m_n=m$. Each of the formal solutions will depend on
a complex parameter $\mu$.

We shall use the notation $h^{(j)}=a_j/\eta-2m_j$. According to our assumptions
we have $h^{(1)}+\ldots+h^{(n)}=0$.

Introduce a new $m$-variable phase function
\begin{equation}\label{phasenew}
\Omega(t_1,\dots,t_m,z_1,\dots,z_n,\tau,p)=
\prod_{j=1}^m\prod_{l=1}^n\Omega_{a_l}(t_j-z_l)
\prod_{1\leq i<j\leq m}\Omega_{-2\eta}(t_i-t_j).
\end{equation}

Introduce a function
\begin{equation}\label{D}
D_M(\mu,z_1,\ldots,z_n)\,=\,\prod_{k=1}^n\,D_k(\mu)^{z_k/p}
\end{equation}
where $D_k$ are defined in \Ref{DD}.

Let $\omega_M(t,\mu,z;p)$ be the function defined in \Ref{ w.f }.
Notice that the parameter $\tau$ of the theta functions in
\Ref{ w.f } is replaced here by $p$ and $\lambda$ is replaced by
$\mu$. It follows from \Ref{b} that
the function $\omega_M(t,\mu,z;p)D_M(\mu,z)$
has the following transformation properties,
\bea 
\kern-20pt&&
\omega_M(t,\mu,\ldots,z_j+p ,\ldots; p)\,D_M(\mu,\ldots,z_j+p,\ldots)
\\[1pt]
\kern-20pt&& =\;e^{2\pi i(\frac{\mu}{2\eta}+m)a_j}\,
\omega_M (t,\mu,\ldots,z_j,\ldots; p)\,D_M(\mu,\ldots,z_j,\ldots)
\eea
for $j=1,\ldots,n$.

For any entire function $\xi$ of one variable, let $E^0_a(z;\xi;M)$,
be the space spanned by the functions of $t\in\C^m$ of the form
\be
\tsize \xi(p\lambda-\Sum_{k=1}^n2a_kz_k+4\eta\Sum_{j=1}^mt_j)
\Omega(t,z)\omega_M(t,\mu,z;p)\,D_M(\mu,z)f(t,z),
\ee
where $f(t,z)$, viewed as a function of $t=(t_1,\dots,t_m)$ belongs to
$\tilde F^m_a(z,\lambda)$ (see \ref{sssym}) for some $\lambda$.

Let $E_a(z;\xi;M)$, the space of cocycles, be the space spanned by functions
of the form $g(t+p\alpha)$, where $g\in E^0_a(z;\xi;M)$ and $\alpha\in\Z^m$.
By construction, $E_a(z;\xi;M)$ is invariant under translations of
the arguments $t_i$ by $p$. We define the space of coboundaries $DE_a(z;\xi)$
to be the subspace of $E_a(z;\xi;M)$ spanned by functions of the form
$f(\dots,t_j+p,\dots)-f(\dots,t_j,\dots)$, $f\in E_a(z;\xi)$.

\begin{proposition}\label{pE(z;M)}
$E_a(z_1,\dots,z_n;\xi;M)=E_a(z_1,\dots,z_j+p,\dots,z_n;\xi;M)$
for $j=1,\dots,n$
\end{proposition}

The proof of this Proposition is the same as the proof of Proposition
\ref{pE(z)}, see \cite{FTV}.

\begin{thm}\label{nir}
Let us write the qKZB equations as $\Psi(\dots,z_j+p,\dots)=K_j(z)\Psi(z)$.
Then, for all entire functions $\xi$, the integrand
\begin{equation}\label{Psi}
e^{-\pi i\frac{\mu\lambda}{2\eta}}\,\Omega(t,z)\,\omega_M(t,\mu,z;p)\,
D_M(\mu,z)\psi^\xi(t,z),
\end{equation}
viewed as a function of $t\in\C^m$, belongs to $E_a(z;\xi;M)\otimes V[0]$ for
all $z\in\C^n$. It obeys the equations
\be
\Psi(t,\dots,z_j+p,\dots)=K_j(z)\Psi(t,z)\mod DE_a(z;\xi;M)\otimes V[0],
\qquad j=1,\dots, n.
\ee
In other words, $\Psi(z,t)$ solves the qKZB equations
in the cohomology
\fixedline
$(E_a(z;\xi;M)/DE_a(z;\xi;M))\otimes V[0]$.
\end{thm}

Consider a square matrix $I$ of size dim $V[0]$ with entries
\ifTwelve
\bea 
&\kern-20pt& I_{L,M}\,={}
\\
&\kern-20pt& e^{-\pi i \frac{\mu\lambda}{2\eta}}\,
\xi(\tau\mu+p\lambda-\Sum_{k=1}^n2a_kz_k+4\eta\Sum_{j=1}^mt_j)
\Omega(t,z)\,\omega_L(t,\lambda,z;\tau)\,\omega_M(t,\mu,z;p)\,
e_L\otimes D_M(\mu,z)e_M,
\eea
\else
\bea 
\ \;I_{L,M}\,=\,e^{-\pi i \frac{\mu\lambda}{2\eta}}\,
\xi(\tau\mu+p\lambda-\Sum_{k=1}^n2a_kz_k+4\eta\Sum_{j=1}^mt_j)
\Omega(t,z)\,\omega_L(t,\lambda,z;\tau)\,\omega_M(t,\mu,z;p)\,
e_L\otimes D_M(\mu,z)e_M,\ifTwelve\kern-1em\fi
\eea
\fi
where $L=(l_1,\ldots,l_n)$, $M=(m_1,\ldots,m_n)$, and
$e_L=e_{l_1}\otimes\ldots\otimes e_{l_n}$,
$e_M=e_{m_1}\otimes\ldots\otimes e_{m_n}$ are basis elements of $V[0]$.
The entries of the matrix are functions of $\tau, p, \lambda, \mu, t, z$.
If we ignore the factor $D_M$, the matrix will be invariant under exchange
of $p,\lambda, L$ and $\tau, \mu, M$.
\begin{corollary}
For every $M$ the corresponding column $I_M=(I_{L,M})$ of the matrix is
a formal solution of the qKZB equations.
\end{corollary}
\begin{proof}
Theorem \ref{nir} follows from the proof of Theorem \ref{irs}. In fact,
according to Theorem \ref{irs} we know that $\Phi(t,z)\psi^\xi(t,z,\lambda)$
is a formal solution. Now we can write
\bea
\kern-20pt&&
\Phi(t,z)\psi^\xi(t,z,\lambda)=e^{-\pi i\frac{\mu\lambda}{2\eta}}\,
e^{\pi i\frac{\mu\lambda}{2\eta}}\,\Phi(t,z)\psi^\xi(t,z,\lambda)={}
\\
\kern-20pt&&
{}=e^{-\pi i\frac{\mu\lambda}{2\eta}}\,\Omega(t,z)\psi^\xi(t,z,\lambda)\,
e^{\pi i\frac{\mu\lambda}{2\eta}}\,
e^{\frac{2\pi i}{p}(m\sum_{l=1}^na_lz_l+\sum_{j=1}^m(2\eta-4\eta j)t_j)}={}
\\
\kern-20pt&&
{}=e^{-\pi i \frac{\mu\lambda}{2\eta}}\,\Omega(t,z)\psi^\xi (t,z,\lambda)\,
e^{\pi i\frac{\mu}{2\eta p}(p\lambda-2\sum_{l=1}^na_lz_l+4\eta\sum_{j=1}^mt_j)}
\times{}
\\
\kern-20pt&&
\quad{}\times e^{\pi i\frac{\mu}{2\eta p}
(2\sum_{l=1}^na_lz_l-4\eta\sum_{j=1}^mt_j)}\,
e^{\frac{2\pi i}{p}(m\sum_{l=1}^na_lz_l+\sum_{j=1}^m(2\eta-4\eta j)t_j)}={}
\\
\kern-20pt&&
{}=e^{-\pi i \frac{\mu\lambda}{2\eta}}\,\Omega(t,z)\psi^\xi(t,z,\lambda)\,
e^{\frac{2\pi i}{p}(\frac{\mu}{2\eta}+m)\sum_{l=1}^na_lz_l}\,
e^{\frac{-2\pi i}{p}\sum_{j=1}^m (\mu +4\eta j-2\eta)t_j}\times{}
\\
\kern-20pt&&
\quad{}\times
e^{\pi i\frac{\mu}{2\eta p}(p\lambda-2\sum_{l=1}^na_lz_l+4\eta\sum_{j=1}^mt_j)}
\,.
\eea
The last factor can be included into the entire function $\xi$.
Thus we proved that
\begin{equation}\label{iis}
e^{-\pi i\frac{\mu\lambda}{2\eta}}\,\Omega(t,z)\psi^\xi(t,z,\lambda)
\,e^{\frac{2\pi i}{p}(\frac{\mu}{2\eta}+m)\sum_{l=1}^na_lz_l}\,
e^{\frac{-2\pi i}{p}\sum_{j=1}^m(\mu +4\eta j-2\eta)t_j}
\end{equation}
is a formal solution depending on an additional parameter $\mu$.

Now choose $M=(m_1,\ldots,m_n)$ as in Theorem \ref{nir}.
The function $\omega_M(t,\mu,z;p)\,D_M(\mu,z)$ has the same transformation
properties under the shifts $t_j\mapsto t_j+p$ and $z_l\mapsto z_l+p$ as
the product of the last two factors in \Ref{iis}. Hence
\begin{equation} 
e^{-\pi i \frac{\mu\lambda}{2\eta}}\,
\Omega(t,z) \psi^\xi (t,z,\lambda)\omega_M(t,\mu,z;p)\,D_M(\mu,z)
\notag
\end{equation}
is a formal integral solution. Theorem \ref{nir} is proved.
\end{proof}

\section{Integration}\label{integration}

\subsection{ Poles of the new formal solutions}

Let $t=(t_1,\ldots,t_m),\,z=(z_1,\ldots,z_n)$. Let $\Om(t,z,\tau,p,\eta)$
be the $m$-variable phase function introduced in \Ref{phasenew}.
Let $\om_{l_1,\dots,l_n}(t,\la,z,\tau,\eta)$ and
$\om_{m_1,\dots, m_n}(t,\mu,z,p,\eta)$ be the functions introduced in
\Ref{ w.f }, here $l_1+\ldots+l_n=m$ and $m_1+\ldots+m_n=m$.

\begin{proposition}\label{Poles}
The poles of the function
\begin{equation}\label{prod}
\Om(t,z,\tau,p,\eta)\,\om_{l_1,\dots, l_n}(t,\la,z,\tau,\eta)\,
\om_{m_1,\dots, m_n}(t,\mu,z,p,\eta)
\end{equation}
are of first order and lie at the hyperplanes given by equations
\bean\label{poles1}
t_i-z_k=-a_k-rp-s\tau+l, \qquad
t_i-z_k=a_k+rp+s\tau+l,
\eean
where $i=1,\ldots,m$, $k=1,\ldots,n$ and $r,s\,\in\Z_{\geq 0}$, $l\in\Z$, and
\bean\label{poles2}
t_i=t_j+2\eta-rp-s\tau+l, \qquad t_i=t_j-2\eta+rp+s\tau+l,
\eean
where $1\leq i<j\leq m$ and $r,s\,\in\Z_{\geq 0}$, $l\in\Z$.
\end{proposition}

\begin{proof}
The poles of the function $\Omega_a(t,\tau, p)$ introduced in \Ref{phase1} are
of first order and are given by
\bean 
t=-a-rp-s\tau+l, \qquad t=a+(r+1)p+(s+1)\tau+l,
\notag
\eean
and the zeros are
\bean 
t=a-rp-s\tau+l, \qquad t=-a+(r+1)p+(s+1)\tau+l,
\notag
\eean
where $r,s\,\in\Z_{\geq 0}$, $l\in \Z$.

Hence the poles of the $m$-variable phase function $\Om(t,z,\tau,p,\eta)$ are
of first order and lie at the hyperplanes defined by equations
\bean 
t_i-z_k=-a_k-rp-s\tau+l, \qquad t_i-z_k=a_k+(r+1)p+(s+1)\tau+l,
\notag
\eean
where $i=1,\ldots,m$, $k=1,\ldots,n$ and $r,s\,\in\Z_{\geq 0}$, $l\in\Z$, and
\bean
t_i=t_j +2\eta -rp-s\tau+l, \qquad t_i=t_j -2\eta +(r+1)p+(s+1)\tau+l,
\notag
\eean
where $1\leq i<j\leq m$ and $r,s\,\in\Z_{\geq 0}$, $l\in\Z$.

The zeros of the $m$-variable phase function $\Om(t,z,\tau,p,\eta)$ lie at
the hyperplanes defined by equations
\bean\label{zeros1}
t_i-z_k=a_k-rp-s\tau+l, \qquad t_i-z_k=-a_k+(r+1)p+(s+1)\tau+l,
\notag
\eean
where $i=1,\ldots,m$, $k=1,\ldots,n$ and $r,s\,\in\Z_{\geq 0}$, $l\in\Z$, and
\bean\label{zeros2}
t_i=t_j -2\eta -rp-s\tau+l, \qquad t_i=t_j +2\eta +(r+1)p+(s+1)\tau+l,
\notag
\eean
where $1\leq i<j\leq m$ and $r,s\,\in\Z_{\geq 0}$, $l\in\Z$.

The poles of the function $\om_{l_1,\dots, l_n}(t,\la,z,\tau,\eta)$ are
of first order and lie at the hyperplanes $t_i-z_k=a_k+ s\tau+l$
where $i=1,\ldots,m$, $k=1,\ldots,n$, $s,l\,\in\Z$, and
at the hyperplanes $t_i=t_j -2\eta+s\tau+l$
where $1\leq i<j\leq m$ and $s,l\in\Z$.

Similarly, the poles of the function $\om_{m_1,\dots, m_n}(t,\mu,z,p,\eta)$
are of first order and lie at the hyperplanes $t_i-z_k=a_k+ sp+l$
where $i=1,\ldots,m$, $k=1,\ldots,n$, $s,l\,\in\Z$, and
at the hyperplanes $ t_i=t_j -2\eta+sp+l$
where $1\leq i<j\leq m$ and $s,l\in\Z$.

Knowing the poles and zeros of the factors in \Ref{prod},
we get the Proposition.
\end{proof}

\subsection{Topology of poles}

The function defined by \Ref{prod} depends on $t,\,z_1,\ldots,z_n$,
$\la,\,\mu$, $\tau,\,p,\,\eta$, $a_1,\ldots,a_n$. Later on we often make
the following assumptions on $\tau,\,p,\,\eta$, $a_1,\ldots,a_n$,
$z_1,\ldots,z_n$.
\begin{equation}\label{im}
\Im\tau>0\,,\qquad \Im p>0\,,\qquad \Im\eta<0\,.
\end{equation}
\begin{equation}\label{indep}
\text{The numbers $\tau$ and $p$ are linearly independent over $\Z$.}
\end{equation}
\begin{equation}\label{eta}
\{2\eta,\, 4\eta,\ldots,\,2m\eta\}\cap\{\Z+\tau \Z+p\Z\}\,=\,\varnothing\,.
\end{equation}
\bean\label{a's}
2a_k+2s\eta\not\in\Z+\tau\Z+p\Z\,,\qquad k=1,\ldots,n\,,\qquad
s=1-m,\ldots,m-1\,.
\eean
\bean\label{z's}
z_l\pm a_l-z_k\pm a_k+2s\eta\not\in\Z+\tau\Z+p\Z\,,\qquad
k,l=1,\ldots,n\,,\quad l\neq k\,,
\\
s=1-m,\ldots,m-1\,,
\notag
\eean
for arbitrary combination of signs.

A set of hyperplanes in an affine space is called a configuration
of hyperplanes. An edge of a configuration is a nonempty intersection
of some hyperplanes of the configuration.

Consider the configuration
${\cal B}={\cal B}(\tau,p,\eta,a_1,\ldots,a_n,z_1,\ldots,z_n)$
of hyperplanes in $\C^m$ defined by equations \Ref{poles1} and \Ref{poles2}.

Fix a natural number $N$ and consider the image of the configuration
${\cal B}$ in $\C^m/N\Z^m$ under the natural projection $\C^m\to\C^m/N\Z^m$.
We shall call the image a configuration of hyperplanes in $\C^m/N\Z^m$ and
denote the image by
${\cal C}_N={\cal C}_N(\tau,p,\eta,a_1,\ldots,a_n,z_1,\ldots,z_n)$.
We always assume that $\Im\tau>0$, $\Im p>0$, and therefore ${\cal C}_N$
is a locally finite collection of hyperplanes in $\C^m/N\Z^m$.

Let $F(t,\la,\mu, z,a,\tau,p)$ be the function defined by \Ref{prod}.
The function is 1-periodic in all variables $t_i$, and therefore it defines
a function on $\C^m/N\Z^m$ which we denote also by $F(t,\la,\mu,z,a,\tau,p)$.
The function $F(t,\la,\mu,z,a,\tau,p)$ is holomorphic on the complement to
the union of hyperplanes of ${\cal C}_N$ in $\C^m/N\Z^m$. Moreover, for any
$\bt\in\Z^m+\tau \Z^m+p\Z^m$, the function $F(t+\bt,\la,\mu,z,a,\tau,p)$ is
holomorphic on the complement to the union of hyperplanes of ${\cal C}_N$
in $\C^m/N\Z^m$.

\begin{proposition}\label{dim}
The number of pairwise distinct edges of any dimension and the dimensions of
all edges of the configuration ${\cal C}_N$ do not depend on the parameters
$\tau, p,$ $\eta, a_1,\ldots,a_n,$ $z_1,\ldots,z_n$ provided that assumptions
\Ref{im}\,--\,\Ref{z's} hold.
\end{proposition}
\begin{proof}
The initial configuration induces a configuration of hyperplanes in any edge
of the initial configuration. The number of pairwise distinct edges of any
dimension and the dimensions of all edges of the initial configuration remain
the same if for each induced configuration the number of its pairwise distinct
hyperplanes does not change. This is obviously true if assumptions
\Ref{im}\,--\,\Ref{z's} hold.
\end{proof}
\begin{corollary}\label{topol}
The topology of the complement of the configuration of hyperplanes ${\cal C}_N$
in $\C^m/N\Z^m$ remains the same for all $\tau,\,p,\,\eta$, $z_1,\ldots,z_n$,
$a_1,\ldots,a_n$ satisfying conditions \Ref{im}\,--\,\Ref{z's}.
\end{corollary}
The Corollary follows from standard reasons in topological singularity theory,
cf.\ \cite{R}.

\subsection{ Hypergeometric integrals }\label{integrals}

Let $\tau, p, \eta, a_1,\ldots,a_n$ be complex numbers.
Assume that \Ref{im} holds and
\begin{equation}\label{Im}
\Im a_k>0\,, \qquad k=1,\ldots,n\,.
\end{equation}

Let $t=(t_1,\ldots,t_m)$, $z=(z_1,\ldots,z_n)$, $a=(a_1,\ldots,a_n)$.
Let $\Om(t,z,a,\tau,p,\eta)$ be the $m$-variable phase function introduced in
\Ref{phasenew}. Let $\om_L(t,\la,z,a,\tau,\eta)$ and $\om_M(t,\mu,z,a,p,\eta)$
be the functions introduced in \Ref{ w.f }, here $L=(l_1,\dots, l_n)$,
$l_1+\ldots+l_n=m$, and $M=(m_1,\dots, m_n)$, $m_1+\ldots+m_n=m$.
Let $\phi(t,\la,\mu,z,a,\tau,p,\eta)$ be a holomorphic function for
$t\in\C^m$, $z,a\in\C^n$, and $\tau,\,p$ lying in the upper half plane.
Assume that $\phi$ is N-periodic in $t$,
$\phi(t+\bt,\la,\mu,z,a,\tau,p,\eta)=\phi(t,\la,\mu,z,a,\tau,p,\eta)$
for all $\bt\in N\Z^m$.

Consider the integral
\bean\label{hyperg}
\kern-20pt&&
I_{LM}^\phi(\la,\mu,z,a,\tau,p,\eta)\,={}
\\[6pt]
\kern-20pt&&
\int_{[0,N]^m}\,\Om(t,z,a,\tau,p,\eta)\,\om_L(t,\la,z,a,\tau,\eta)\,
\om_M(t,\mu,z,a,p,\eta)\,\phi(t,\la,\mu,z,a,\tau,p,\eta)\,dt
\notag
\eean
where $dt=dt_1\wedge\ldots\wedge dt_n$.
The integral will be called {\it a hypergeometric integral}.

The integrand, $F(t,\la,\mu,z,a,\tau,p,\eta)$, in \Ref{hyperg} is N-periodic,
hence the integral can be considered as an integral over the image
in $\C^m/N\Z^m$ of the subspace $\R^m$ under the natural projection.

For any $\bt\in\C^m$, denote by
$I_{LM}^\phi(\la,\mu,z,a,\tau,p,\eta)_\bt$ the integral
\begin{equation}\label{hyperg.al}
\int_{[0,N]^m}\,F(t+\bt,\la,\mu,z,a,\tau,p)\,dt\,.
\end{equation}

Fix $\bt\in\tau\Z^m+p\Z^m$ and $z\in\C^n$. Assume that for all $i=1,\ldots,n$,
\bean \label{hor}
\Im a_i\gg\Im\tau\,,\quad\ \Im a_i\gg\Im p\,,\quad\ \text{and}\ \quad
-\Im\eta\gg\Im\tau\,,\quad\ -\Im\eta\gg\Im p\,.
\eean
Then it is easy to see that the poles of the integrand of \Ref{hyperg.al} lie
far from the integration cycle $[0,N]^m$. Hence the integral is well defined
and holomorphically depends on the parameters $\la,\mu,z,a,\tau,p,\eta$.
We shall call this range
of parameters {\it the starting range for given $\bt$ and $z$.}

\begin{proposition}\label{indepe}
For fixed $\bt,\gm\in\tau\Z^m+p\Z^m$ and $z\in \C^n$, the integrals
$I_{LM}^\phi(\la,\mu,z,a,\tau,p,\eta)_\bt$ and
$I_{LM}^\phi (\la,\mu,z,a,\tau,p,\eta)_\gm$ are equal
if condition \Ref{hor} hold.
\end{proposition}

\begin{proof}
Change variables, $t\to t'-\bt+\gm$, in the first integral. Then the integrand
of the first integral becomes equal to the integrand of the second while the
integration cycle becomes equal to $[0,N]^m +\bt-\gm$. If condition \Ref{hor}
hold, then the integrand has no poles around the tori $[0,N]^m+\bt-\gm$ and
$[0,N]^m$. Since the integrand is a closed differential form, the integrals
over $[0,N]^m +\bt-\gm$ and $[0,N]^m$ are equal.
\end{proof}

So far we defined the function $I_{LM}^\phi(\la,\mu,z,a,\tau,p,\eta)_\bt$ in
the starting range of parameters with respect to given $\bt$ and $z$ as
the integral of ${F(t+\bt,\la,\mu,z,a,\tau,p)\,dt}$ over the torus
$[0,N]^m\subset\C^m/N\Z^m$. Moreover, this function does not depend on $\bt$.
In order to define the function $I_{LM}^\phi(\la,\mu,z,a,\tau,p,\eta)_\bt$
for all values of parameters satisfying conditions \Ref{im}\,--\,\Ref{z's}
we use analytic continuation. The result of the analytic continuation can be
represented as an integral of the integrand over a suitably deformed torus,
which we denote by $T^m_N$. Namely, the poles of the integrand are located
at the hyperplanes of the configuration
${\cal C}_N(\tau,p,\eta,a_1,\ldots,a_n,z_1,\ldots,z_n)$. We deform
the parameters $\tau,\,p,\,\eta$, $a_1,\ldots,a_n$, $z_1,\ldots,z_n$ preserving
conditions \Ref{im}\,--\,\Ref{z's}. Then by Corollary \ref{topol}, the topology
of the complement of ${\cal C}_N$ in $\C^m/N\Z^m$ does not change. We deform
the integration torus accordingly so that the torus does not intersect
the hyperplanes of the configuration at every moment of the deformation.
Then the analytic continuation of the function
$I_{LM}^\phi(\la,\mu,z,a,\tau, p,\eta)_\bt$ is given by the integral
\begin{equation}\label{Hyperg.al}
\int_{T^m_N}\,F(t+\bt,\la,\mu,z, a,\tau,p)\,dt\,.
\end{equation}

\begin{thm}\label{analytic}
The integral \Ref{hyperg.al} can be analytically continued as a holomorphic
univalued function to the domain of the parameters $\la,\mu,z,a,\tau,p,\eta$
satisfying conditions \Ref{im}\,--\,\Ref{z's}.
\end{thm}
The proof of the Theorem is the same as the proof of Theorem 5.7 in \cite{TV1}.

We have the following important corollary.

\begin{corollary}\label{indepen}
The function $I_{LM}^\phi(\la,\mu,z,a,\tau,p,\eta)_\bt$
defined by \Ref{Hyperg.al} does not depend on $\bt$.
\end{corollary}

\subsection{Solutions to qKZB}\label{sol.qKZB}
Fix complex numbers $\tau,\,p,\,\eta$, $\Lambda_1,\ldots,\Lambda_n$ and set
$a_i=\eta\Lambda_i$. Assume that the parameters $\tau,\,p,\,\eta$,
$a_1,\ldots,a_n$ satisfy conditions \Ref{im}\,--\,\Ref{a's} and
$\Lambda_1+\ldots+\Lambda_n=m$ for some positive integer $m$.

Let $\Om(t,z,a,\tau,p,\eta)$ be the $m$-variable phase function introduced in
\Ref{phasenew}. Let $\om_L(t,\la,z,a,\tau,\eta)$ and $\om_M(t,\mu,z,a,p,\eta)$
be the functions introduced in \Ref{ w.f }, here $L=(l_1,\dots,l_n)$,
$|L|=l_1+\ldots+l_n=m$, and $M=(m_1,\dots, m_n)$, $|M|=m$.

Fix a natural number $N$. Let $\xi$ be an entire function of one variable
which is $4\eta N$-periodic, $\xi(\la+4\eta N)=\xi(\la)$.

Let $V_\Lambda=\oplus_{j=0}^\infty\C e_j$ and
$V=V_{\Lambda_1}\otimes\ldots\otimes V_{\Lambda_n}$.
For any $L=(l_1,\dots,l_n)$, $|L|= m$, set
$e_L=e_{l_1}\otimes\ldots\otimes e_{l_n}{{}\,\in\,V}$.

Introduce a $V[0]\otimes V[0]$-valued function $u^\xi$ by
\bean\label{sol}
u^\xi(z_1,\ldots,z_n,\la,\mu,\tau,p)\,=\,\sum_{L,M,\,|L|=|M|=m}\,
u^\xi_{LM}(z_1,\ldots,z_n,\la, \mu, \tau,p)\, e_L\otimes e_M
\eean
where
\bean\label{coord}
\kern-20pt&&
u^\xi_{LM}(z_1,\ldots,z_n,\la,\mu,\tau,p)\,=\,
e^{-\pi i{\mu\la\over2\eta}}\,\times
\\[6pt]
\kern-20pt&&
\int_{T^m_N}\,\xi(p\la+\tau\mu-
\tsize\Sum_{l=1}^n 2a_lz_l+4\eta\Sum_{j=1}^m t_j)\,
\Om(t,z,a,\tau,p,\eta)\,\om_L(t,\la,z,a,\tau,\eta)\,
\om_M(t,\mu,z,a,p,\eta)\,dt\,.\ifTwelve\kern-1em\fi
\notag
\eean
Here we assume that $z_1,\ldots,z_n$ satisfy condition \Ref{z's} and we define
the integral by analytic continuation described in Section \ref{integrals}.

For any $j=1,\ldots,n$, introduce an $\End(V[0])$-valued function $D_j(\la)$
by the formula
\bean\label{D.op}
D_j(\la):e_L\,\mapsto\,\,
\frac{\al(\la\,-\,2\eta\,\sum_{l=1}^{j}h^{(l)})}
{\al(\la\,-\,2\eta\,\sum_{l=1}^{j-1}h^{(l)})}\,
e^{\pi i\eta\Lambda_j(\sum_{l=1}^{j-1}\Lambda_l-\sum_{l=j+1}^{n}\Lambda_l)}
\,e_L\,,
\eean
cf.\ \Ref{DD}, and an $\End(V[0])$-valued function $D(\la,z,p)$ by
\begin{equation}\label{De}
D(\la,z,p)\,=\,\prod_{j=1}^n\,D_{j}(\mu)^{z_j/p}\,,
\end{equation}
cf.\ \Ref{D}.

\begin{thm}\label{u-solu}
Under the above conditions, for any $j=1,\ldots,n$, we have
\bean\label{eq1}
u(\ldots,z_j+p,\ldots)=K_j(z_1,\ldots,z_n,\la,\tau,p)
\otimes D^{-1}_j(\mu)\,\,u(\ldots,z_j,\ldots),
\eean
\bean\label{eq2}
u(\ldots,z_j+\tau,\ldots)=D_j^{-1}(\la)\otimes K_j(z_1,\ldots,z_n,\mu,p,\tau)
\,\,u(\ldots,z_j,\ldots),
\eean
and, if in addition the function $\xi$ is $2a_l$-periodic
for all $l=1,\ldots,n$, then
\bean\label{eq3}
u(\ldots,z_j+1,\ldots)=u(\ldots,z_j,\ldots).
\eean
Here $K_j(z_1, \ldots,z_n,\la,\tau, p)$ is the qKZB operator defined
by \Ref{qKZB}, i.e.\ the operator of the qKZB equations with step $p$ and
defined in terms of elliptic $R$-matrices with modulus $\tau$.
\end{thm}

The following Corollary is equivalent to the first two statements of
the Theorem. The third statement is trivial.

\begin{corollary}\label{solu}
Let $f:V\to \C$ be a linear function.
Consider the functions
\bean\label{solu1}
\Psi_f^\xi(z,\la,\mu,\tau,p)\,=\,
(1\otimes f)\,(1\otimes D(\mu,z,p))\,u(z,\la,\mu,\tau,p)\,,
\eean
and
\bean\label{solu2}
\Phi_f^\xi(z,\la,\mu,\tau,p)\,=\,
(f\otimes 1)\,(D(\la,z,\tau)\otimes 1)\,u(z,\la,\mu,\tau,p)\,.
\eean
Then for a fixed $\mu$, the function $\Psi^\xi$ is a solution
of the qKZB equations with modulus $\tau$ and step $p$,
\bean\label{solu3}
\kern-20pt&&
\Psi_f^\xi(z_1,\ldots,z_j+p,\ldots,z_n,\la,\mu,\tau,p)\,={}
\\
\kern-20pt&&
{}=\,K_j(z_1,\ldots,z_n,\la,\tau, p)\,
\Psi_f^\xi(z_1,\ldots,z_n,\la,\mu,\tau,p)\,,\qquad j=1,\ldots,n\,,
\notag
\eean
and for a fixed $\la$, the function $\Phi^\xi$ is a solution
of the qKZB equations with modulus $p$ and step $\tau$,
\bean\label{solu4}
\kern-20pt&&
\Phi_f^\xi(z_1,\ldots,z_j+\tau,\ldots,z_n,\la,\mu,\tau,p)\,={}
\\
\kern-20pt&&
{}=K_j(z_1, \ldots,z_n,\mu,p,\tau)\,
\Psi_f^\xi(z_1,\ldots,z_n,\la,\mu,\tau,p)\,,\qquad j=1,\ldots,n\,.
\notag
\eean
\end{corollary}

\begin{proof}
Fix $M=(m_1,\dots,m_n)$, $|M|=m$. By Theorem \ref{nir}, the function $\Psi$
defined by \Ref{Psi} is a formal solution to the qKZB equations \Ref{solu3}.
This means that for any $j$, the difference
$$
\Psi^\xi(z_1,\ldots,z_j+p,\ldots,z_n,\la,\mu,\tau,p)-
K_j(z_1,\ldots,z_n,\la,\tau, p)\,
\Psi^\xi(z_1,\ldots,z_n,\la,\mu,\tau,p)
$$
is an element of the space $DE_a(z;\xi;M)\otimes V[0]$. By Corollary
\ref{indepen}, the integral in \Ref{coord}, defined by the analytic
continuation described in Section \ref{integrals}, is equal to zero on
$DE_a(z;\xi;M)\otimes V[0]$. Hence, if $f:V\to \C$ is a linear function,
then, for a fixed $\mu$, the function $\Psi_f^\xi$, defined by \Ref{solu1},
is a solution to the qKZB equations \Ref{solu3}. This proves the first
statement of Corollary \ref{solu}. The proof of the second statement is
similar. The Corollary implies the Theorem.
\end{proof}

The solutions of the qKZB equations constructed in Corollary \ref{solu} depend
on an entire function $\xi$ which is $4\eta N$-periodic. An important example
of such a function is the function $\xi \equiv 1$. In the next section we
consider another important example of entire functions.

\subsection{Theta function properties of solutions}\label{Theta}
Recall that a scalar theta function of level $N$ is a function $f$ such that
$f(\la +1)=f(\la)$, \ $f(\la+\tau)=e^{-\pi iN(2\la+\tau)}f(\la)$.

Consider the objects described in Section \ref{sol.qKZB}, in particular,
let $V=V_{\Lambda_1}\otimes\ldots\otimes V_{\Lambda_n}$ as in
Section \ref{sol.qKZB}. We say that a $V[0]$-valued function $\Psi(\la)$
is a vector-valued theta function of level $N$, if
\bean\label{theta1}
\Psi(\la+1)\,=\,\Psi(\la),
\eean
\bean\label{theta2}
\Psi(\la+\tau)\,=\,e^{-\pi i(2 N\la+N\tau+\sum_{j=1}^n h^{(j)}
(z_j-a_j-2a_{j+1}-\ldots-2a_{n}))}\,
\Psi(\la).
\eean
Clearly, the space of holomorphic vector-valued theta functions of level $N$
is finite-dimensional.

Assume that the numbers $p$ and $\eta$ are such that $-p/4\eta=N$ where $N$
is a natural number. Let $f(\la)$ be a scalar theta function of level $N$.
Set $\xi(\la)=f(\la/p)$. Then $\xi$ is an entire $4\eta N$-periodic function.

Consider the solution $\Psi^\xi_f(z,\la,\mu,\tau,p,\eta)$ to the qKZB equations
constructed in Corollary \ref{solu}.

\begin{thm}\label{theta.sol}
Under the above assumptions, if the parameter $\mu$ has the form
$\mu=2\eta(m+2s)$, $s\in\Z$, then the solution
$\Psi^\xi_f(z,\la,\mu,\tau,p,\eta)$ is a vector-valued theta function of
level $m+N$ as a function of variable $\la$.
\end{thm}
\begin{proof}
We have $\Psi^\xi_f(\ldots,\la+1,\ldots)=e^{-\pi i\mu/2\eta}(-1)^m
\Psi^\xi_f(\ldots,\la,\ldots)$ where the first factor comes from the factor
$e^{-\pi i{\mu\la\over2\eta}}$ in \Ref{coord} and the second factor comes from
$\om_L(\ldots,\la,\ldots)$ in \Ref{coord}. Hence, according to our assumptions,
$\Psi^\xi_f(\ldots,\la,\ldots)$ satisfies \Ref{theta1}.

In order to check property \Ref{theta2} we shall check the property for each
coordinate function $u^\xi_{LM}$ defined by \Ref{coord}.

We have $u^\xi_{LM}(\ldots,\la+\tau,\ldots)=e^{-\pi iA}e^{-\pi iB}e^{-\pi iC}
u^\xi_{LM}(\ldots,\la,\ldots)$ where the first, second and third factors
correspond to transformations of $e^{-\pi i\mu/2\eta}$, $\xi$ and $\om_L$ in
\Ref{coord}, respectively. Here $A=\mu\tau/2\eta$,
$B=N(2\la+2\tau\mu/p-4\eta/p\,\sum\Lambda_kz_k+8\eta/p\,\sum t_j+\tau)$,
$$
C\,=\,m\tau\,+\,\sum_{k=1}^n \,\sum_{j=l_1+\ldots+l_{k-1}+1}^{l_1+\ldots+l_k}
2\bigl(\la+t_j-z_k-a_k+2l_k\eta-
\tsize 2\eta\Sum_{q=1}^{k-1}(\La_q-2l_q)\bigr)\,.
$$
Using our
assumptions, after simple calculations we get
$$
A+B+C=2(N+m)\la+(N+m)\tau+
\tsize\Sum_{j=1}^n h^{(j)} (z_j-a_j-2a_{j+1}-\ldots-2a_{n})
$$
where $h^{(j)}=\La_j-2l_j$. The Theorem is proved.
\end{proof}

\section{Monodromy of solutions of qKZB equations}\label{Monodromy}

\subsection{ Monodromy with respect to permutations of variables}

\begin{thm}\label{Perm}
Let $\Psi (z_1,\ldots,z_n,\la)$ be a solution of the qKZB equations with values
in $(V_{\La_1}\otimes \ldots\otimes V_{\La_n})[0]$ step $p$ and modulus $\tau$.
Then for any $j=1,\ldots,n-1$, the function
\ifTwelve
\bean\label{perm}
\kern-20pt&&
\Psi_{j}(z_1,\ldots,z_n,\la)\,={}
\\
\kern-20pt&&
{}=\,P^{(j,j+1)}\,R_{\La_j, \La_{j+1}}^{(j,j+1)}
(z_{j+1}-z_j,\la-2\eta\tsize\Sum_{l=1}^{j-1} h^{(l)},\tau)\,
\Psi(z_1,\ldots,z_{j+1},z_j,\ldots ,z_n,\la)
\notag
\eean
\else
\bean\label{perm}
\Psi_{j}(z_1,\ldots,z_n,\la)\,=\,P^{(j,j+1)}\,R_{\La_j, \La_{j+1}}^{(j,j+1)}
(z_{j+1}-z_j,\la-2\eta\tsize\Sum_{l=1}^{j-1} h^{(l)},\tau)\,
\Psi(z_1,\ldots,z_{j+1},z_j,\ldots ,z_n,\la)
\eean
\fi
is a solution of the qKZB equations with values in $(V_{\La_1}\otimes\ldots
V_{\La_{j+1}}\otimes V_{\La_j}\ldots\otimes V_{\La_n})[0]$, step $p$ and
modulus $\tau$. Here $P^{(j,j+1)}$ is the permutation of the $j$-th and
$j+1$-th factors, and $R_{\La_j,\La_{j+1}}(z,\la,\tau)\in
\End(V_{\La_{j}}\otimes V_{\La_{j+1}})$ is the elliptic $R$-matrix
with modulus $\tau$.
\end{thm}
\begin{proof}
The proof of the Theorem is straightforward. Let us check, for instance,
that $\Psi_j$ satisfies the first qKZB equation assuming that $j>1$.
We have
\bean
\kern-20pt&&
\Psi_j(z_1+p,\ldots)\,=\,P^{(j,j+1)}\,R_{\La_j,\La_{j+1}}^{(j,j+1)}
(z_{j+1}-z_j,\la-2\eta\tsize\Sum_{l=1}^{j-1} h^{(l)},\tau)\,
\notag
\\
\kern-20pt&&
{}\times K_1^{\La_1,\ldots,\La_n}(z_1,\ldots,z_{j+1},z_j,\ldots,z_n,\la)\,
\Psi(z_1,\ldots,z_{j+1},z_j,\ldots )\,={}
\notag
\\
\kern-20pt&&
=\,\,P^{(j,j+1)}\,R_{\La_j,\La_{j+1}}^{(j,j+1)}
(z_{j+1}-z_j,\la-2\eta\tsize\Sum_{l=1}^{j-1} h^{(l)},\tau)\,\Gm_1\,
R_{\La_1,\La_{n}}^{(1,n)}(z_1-z_n,\la-2\eta\Sum_{l=2}^{n-1} h^{(l)},\tau)
\,\ldots{}
\notag
\\
\kern-20pt&&
{}\times R_{\La_1, \La_{j+1}}^{(1,j+1)}
(z_1-z_j,\la-2\eta\tsize\Sum_{l=2}^{j} h^{(l)},\tau)\,
R_{\La_1, \La_{j}}^{(1,j)}(z_1-z_{j+1},\la-2\eta\Sum_{l=2}^{j-1} h^{(l)},\tau)
\,\ldots{}
\notag
\\
\kern-20pt&&
{}\times R_{\La_1, \La_{2}}^{(1,2)}(z_1-z_2,\la,\tau)\,
\Psi(z_1,\ldots,z_{j+1},z_j,\ldots)\,={}
\notag
\\
\kern-20pt&&
=\,\,P^{(j,j+1)}\,\Gm_1\,R_{\La_1,\La_{n}}^{(1,n)}
(z_1-z_n,\la-2\eta\tsize\Sum_{l=2}^{n-1} h^{(l)},\tau)\,\ldots\,
R_{\La_j,\La_{j+1}}^{(j,j+1)}
(z_{j+1}-z_j,\la-2\eta\Sum_{l=2}^{j-1} h^{(l)},\tau)
\notag
\\
\kern-20pt&&
{}\times R_{\La_1,\La_{j+1}}^{(1,j+1)}
(z_1-z_j,\la-2\eta\tsize\Sum_{l=2}^{j} h^{(l)},\tau)\,
R_{\La_1, \La_{j}}^{(1,j)}(z_1-z_{j+1},\la-2\eta\Sum_{l=2}^{j-1} h^{(l)},\tau)
\,\ldots{}
\notag
\\
\kern-20pt&&
{}\times R_{\La_1,\La_{2}}^{(1,2)}(z_1-z_2,\la,\tau)\,
\Psi(z_1,\ldots,z_{j+1},z_j,\ldots)\,={}
\notag
\\
\noalign{\goodbreak}
\kern-20pt&&
=\,\,P^{(j,j+1)}\,\Gm_1\,R_{\La_1,\La_{n}}^{(1,n)}
(z_1-z_n,\la-2\eta\tsize\Sum_{l=2}^{n-1} h^{(l)},\tau)\,\ldots\,
R_{\La_1, \La_{j}}^{(1,j)}
(z_1-z_{j+1},\la-2\eta\Sum_{l=2,\,l\neq j}^{j+1} h^{(l)},\tau)\,\ldots{}
\notag
\\
\kern-20pt&&
{}\times R_{\La_1,\La_{j+1}}^{(1,j+1)}
(z_1-z_j,\la-2\eta\tsize\Sum_{l=2}^{j-2} h^{(l)},\tau)\,
R_{\La_j,\La_{j+1}}^{(j,j+1)}
(z_{j+1}-z_j,\la-2\eta\Sum_{l=1}^{j-1} h^{(l)},\tau)
\notag
\\
\kern-20pt&&
{}\times R_{\La_1,\La_{2}}^{(1,2)}(z_1-z_2,\la,\tau)\,
\Psi(z_1,\ldots,z_{j+1},z_j,\ldots )\,={}
\notag
\\[4pt]
\kern-20pt&&
=\,\,K_1^{\La_1,\ldots,\La_{j+1},\La_j,\ldots,\La_n}
(z_1,\ldots,z_j,z_{j+1},\ldots,z_n,\la)\,\Psi_j (z_1,\ldots )\,.
\notag
\eean
The other cases are checked similarly.
\end{proof}

\noindent{\sl Remark.} The transformation
\ifTwelve
\bea 
\kern-20pt&& T_j:\Psi(z_1,\ldots,z_n,\la)\to
\\[2pt]
\kern-20pt&&\quad P^{(j,j+1)}\,R_{\La_j, \La_{j+1}}^{(j,j+1)}
(z_{j+1}-z_j,\la-2\eta\tsize\Sum_{l=1}^{j-1} h^{(l)},\tau)\,
\Psi(z_1,\ldots,z_{j+1},z_j,\ldots ,z_n,\la)
\eea
\else
\bea 
T_j:\Psi(z_1,\ldots,z_n,\la)\to P^{(j,j+1)}\,R_{\La_j, \La_{j+1}}^{(j,j+1)}
(z_{j+1}-z_j,\la-2\eta\tsize\Sum_{l=1}^{j-1} h^{(l)},\tau)\,
\Psi(z_1,\ldots,z_{j+1},z_j,\ldots ,z_n,\la)
\eea
\fi
preserves the theta function properties \Ref{theta1}, \Ref{theta2}.

\begin{proposition}\label{theta.perm}
Assume that a $(V_{\La_1}\otimes \ldots \otimes V_{\La_n})[0]$-valued function
$\Psi(z_1,\ldots,z_n,\la)$ satisfies the theta function properties
\Ref{theta1}, \Ref{theta2}. Then for any $j=1,\ldots,n$,
the $(V_{\La_1}\otimes\ldots V_{\La_{j+1}}\otimes V_{\La_j}\ldots\otimes
V_{\La_n})[0]$-valued function $T_j\Psi$ also satisfies the theta function
properties \Ref{theta1}, \Ref{theta2}.
\end{proposition}

The Proposition easily follows from Proposition \ref{l-trans}.

\vs
Let $u^{\xi,\La_1,\ldots,\La_n}(z_1,\ldots,z_n,\la,\mu,\tau,p)$ be
the function constructed for the tensor product
$V=V_{\La_1}\otimes\ldots\otimes V_{\La_n}$ in \Ref{sol}.
Let $D_{\La_1,\ldots,\La_n}(\mu,z,p)$ be the $\End(V[0])$-valued function
defined in \Ref{De}. According to Theorem \ref{u-solu},
the $V[0]\otimes V[0]$-valued function
\bean\label{resh}
\Psi^{\La_1,\ldots,\La_n}(\mu,z)\,=\,
(1\otimes D_{\La_1,\ldots,\La_n}(\mu, z,p))\,
u^{\xi,\La_1,\ldots,\La_n}(z_1,\ldots,z_n,\la,\mu,\tau,p)
\eean
is a solution of the qKZB equations with respect to the first factor,
\ifTwelve
\bea
\kern-20pt&& \Psi^{\La_1,\ldots,\La_n}(\ldots,z_k+p,\ldots)\,={}
\\
\kern-20pt&& (K_k(z_1,\ldots,z_n,\la,\tau,p)\otimes 1)\,\,
\Psi^{\La_1,\ldots,\La_n}(\ldots,z_k,\ldots)\,,\qquad k=1,\ldots,n.
\eea
\else
\bea
\quad\Psi^{\La_1,\ldots,\La_n}(\ldots,z_k+p,\ldots)\,=\,
(K_k(z_1,\ldots,z_n,\la,\tau,p)\otimes 1)\,\,
\Psi^{\La_1,\ldots,\La_n}(\ldots,z_k,\ldots)\,,\qquad k=1,\ldots,n.
\eea
\fi

For any $j=1,\ldots,n$, denote by $V_j$ the tensor product
$V_{\La_1}\otimes\ldots V_{\La_{j+1}}\otimes V_{\La_j}\ldots\otimes V_{\La_n}$.
According to Theorem \ref{Perm}, for any $j$, the $V[0]\otimes V_j[0]$-valued
function
\bean\label{u-perm}
\kern-20pt&&
\Psi_j(z,\la,\mu,\tau,p)\,=\,
(P^{(j,j+1)}R_{\La_{j+1}, \La_{j}}^{(j,j+1)}
(z_{j+1}-z_j,\la-2\eta\tsize\Sum_{l=1}^{j-1} h^{(l)},\tau)\otimes{}
\notag
\\
\kern-20pt&&
D_{\ldots,\La_{j+1},\La_j,\ldots}(\mu,\ldots,z_{j+1},z_j,\ldots,p ))\,
u^{\xi,\ldots,\La_{j+1},\La_j,\ldots}(\ldots,z_{j+1},z_j,\ldots)
\notag
\eean
is a solution with respect to the first factor of the same qKZB equations.

The next Theorem describes a relation between the two solutions and can be
considered as a description of the monodromy of the hypergeometric solutions
constructed in Section \ref{sol.qKZB} with respect to permutation of variables.

Introduce a new $R$-matrix $\tilde R_{A,B}(z,\mu,p)\in\End(V_A\otimes V_B)$ by
\ifTwelve
\bean\label{new-r}
\kern-20pt&&
\tilde R_{A,B} (z,\mu,p)\,={}
\\
\kern-20pt&&
e^{2\pi i ABz/p}\,\Bigl(\,\frac{\al(\la)}{\al(\la-2\eta h^{(2)})}\Bigr)^{z/p}\,
R_{A,B} (z,\mu, p)\,\Bigl(\,\frac{\al(\la-2\eta(h^{(1)}+h^{(2)})}
{\al(\la -2\eta h^{(1)})}\Bigr)^{z/p}.
\notag
\eean
\else
\bean\label{new-r}
\tilde R_{A,B} (z,\mu,p)\,=\,
e^{2\pi i ABz/p}\,\Bigl(\,\frac{\al(\la)}{\al(\la-2\eta h^{(2)})}\Bigr)^{z/p}\,
R_{A,B} (z,\mu, p)\,\Bigl(\,\frac{\al(\la-2\eta(h^{(1)}+h^{(2)})}
{\al(\la -2\eta h^{(1)})}\Bigr)^{z/p}.
\eean
\fi

\begin{thm}\label{permut}
\ifTwelve
\bean\label{perm-1}
\kern-20pt&&
\Psi_j(z,\la,\mu,\tau,p)\,={}
\\
\kern-20pt&&
(1\otimes P^{(j,j+1)}
\tilde R_{\La_{j}, \La_{j+1}}^{(j,j+1)}
( z_{j}- z_{j+1}, \mu - 2\eta \tsize\Sum_{l=1}^{j-1} h^{(l)},p))
\Psi^{\La_1,\ldots,\La_n}(z,\la,\mu,\tau,p)\,.
\notag
\eean
\else
\bean\label{perm-1}
\Psi_j(z,\la,\mu,\tau,p)\,=\,(1\otimes P^{(j,j+1)}
\tilde R_{\La_{j}, \La_{j+1}}^{(j,j+1)}
( z_{j}- z_{j+1}, \mu - 2\eta \tsize\Sum_{l=1}^{j-1} h^{(l)},p))
\Psi^{\La_1,\ldots,\La_n}(z,\la,\mu,\tau,p)\,.
\eean
\fi
\end{thm}

\noindent{\sl Remark.} According to Proposition \ref{transR}, the matrix
$\tilde R_{A,B}(z,\mu,p)$ is $p$-periodic,
$$
\tilde R_{A,B}(z+p,\mu,p)\,=\,\tilde R_{A,B}(z,\mu,p)\,.
$$
Hence, formula \Ref{perm-1} expresses solution $\Psi_j$ as a linear combination
of solutions $\Psi^{\La_1,\ldots,\La_n}$ with $p$-periodic coefficients.

\vs
\begin{proof}
First let us prove that
\bean\label{Permut}
\kern-30pt&&
\Psi_j(z,\la,\mu,\tau,p)\,=\,\bigl(1\otimes
D_{\ldots,\La_{j+1},\La_j,\ldots}(\mu,\ldots,z_{j+1},z_j,\ldots,p)
\\
\kern-30pt&&
{}\times P^{(j,j+1)}\,R_{\La_{j},\La_{j+1}}^{(j,j+1)}
(z_{j}-z_{j+1},\mu-2\eta \tsize\Sum_{l=1}^{j-1} h^{(l)},p)
\notag
\\
\kern-30pt&&
{}\times
D_{\ldots,\La_{j},\La_{j+1},\ldots}^{-1}(\mu,\ldots,z_{j},z_{j+1},\ldots,p))
\,\Psi^{\La_1,\ldots,\La_n}(z,\la,\mu,\tau,p)\,.
\notag
\eean
In fact, we have
\bean\label{u}
\kern-30pt&&
\Psi_j(z,\la,\mu,\tau,p)\,=\,
(P^{(j,j+1)}R_{\La_{j+1}, \La_{j}}^{(j,j+1)}
(z_{j+1}-z_j,\la-2\eta \tsize\Sum_{l=1}^{j-1} h^{(l)},\tau)\otimes{}
\\
\kern-30pt&&
D_{\ldots,\La_{j+1},\La_j,\ldots}(\mu,\ldots,z_{j+1},z_j,\ldots))\,
\,\sum_{L,M,\,|L|=|M|=m}u_{LM}(\ldots,z_{j+1},z_j,\ldots)\,e_L\otimes e_M
\notag
\eean
where $L=(l_1,\ldots,l_{j+1},l_j,\ldots,l_n)$,
$M=(m_1,\ldots,m_{j+1},m_j,\ldots,m_n)$, and
\bea 
u_{LM}(\ldots,z_{j+1},z_j,\ldots)\,=\,e^{-\pi i{\mu\la\over2\eta}}\,
\int_{T^m_N}\,\xi(p\la+\tau\mu-\tsize\Sum_{l=1}^n 2a_lz_l+
4\eta\Sum_{j=1}^m t_j)\,\Om(t,z,a,\tau,p,\eta)
\\
{}\times\,
\om^{\,\ldots,\La_{j+1},\La_{j},\ldots}_L(t,\la,\ldots,z_{j+1},z_j,\ldots,\tau)
\,\om^{\,\ldots,\La_{j+1},\La_{j},\ldots}_M(t,\mu,\ldots,z_{j+1},z_j,\ldots,p)
\,dt\,.\!\!\!
\eea
According to \Ref{R.c}, we have
\ifTwelve
\bean\label{a}
\kern-20pt&& \!\!P^{(j,j+1)}\,\sum_L\,\om^{\ldots,\La_{j+1},\La_{j},\ldots}
_L(t,\la, \ldots,z_{j+1},z_j,\ldots,\tau)
\\
\kern-20pt&& {}\times\,R_{\La_{j+1}, \La_{j}}^{(j,j+1)}
(z_{j+1}-z_j,\la-2\eta{\tsize\Sum_{l=1}^{j-1}} h^{(l)},\tau)\,e_L\;
=\,\sum_S\,\om^{\ldots,\La_{j},\La_{j+1},\ldots}
_S(t,\la,\ldots,z_{j},z_{j+1},\ldots,\tau)\,e_S
\notag
\eean
\else
\bean\label{a}
P^{(j,j+1)}\,\sum_L\,\om^{\ldots,\La_{j+1},\La_{j},\ldots}
_L(t,\la, \ldots,z_{j+1},z_j,\ldots,\tau)\,R_{\La_{j+1}, \La_{j}}^{(j,j+1)}
(z_{j+1}-z_j,\la-2\eta\tsize\Sum_{l=1}^{j-1} h^{(l)},\tau)\,e_L
\\
=\,\sum_S\,\om^{\ldots,\La_{j},\La_{j+1},\ldots}
_S(t,\la,\ldots,z_{j},z_{j+1},\ldots,\tau)\,e_S
\notag
\eean
\fi
where $e_S\in V$ and
\bean\label{g}
\kern-16pt&&\kern-.5em \sum_M\,\om^{\ldots,\La_{j+1},\La_{j},\ldots}
_M(t,\mu, \ldots,z_{j+1},z_{j},\ldots ,p)\,e_M\,=
\\
\kern-16pt&& =\,P^{(j,j+1)}\,\sum_Q\,\om^{\ldots,\La_{j},\La_{j+1},\ldots}
_Q(t,\mu, \ldots,z_{j},z_{j+1},\ldots,p)\,R_{\La_{j},\La_{j+1}}^{(j,j+1)}
(z_{j}-z_{j+1},\mu-2\eta\tsize\Sum_{l=1}^{j-1} h^{(l)},p)\,e_Q\,
\ifTwelve\kern-1em\fi
\notag
\eean
where $e_Q\in V$. Formulae \Ref{a}, \Ref{g} and \Ref{u} prove \Ref{Permut}.

Now the Theorem follows from the following Lemma.

\begin{lemma}
\bea
\kern-20pt&&
D_{\ldots,\La_{j+1},\La_j,\ldots}(\mu,\ldots,z_{j+1},z_j,\ldots,p)
P^{(j,j+1)}R_{\La_{j}, \La_{j+1}}^{(j,j+1)}
(z_{j}-z_{j+1},\mu-2\eta\Sum_{l=1}^{j-1} h^{(l)},p)
\\
\kern-20pt&&
\times\,
D_{\ldots,\La_{j},\La_{j+1},\ldots}^{-1}(\mu,\ldots,z_{j},z_{j+1},\ldots,p ))\,
=\,P^{(j,j+1)}\tilde R_{\La_{j},\La_{j+1}}^{(j,j+1)}
(z_{j}-z_{j+1}, \mu-2\eta\Sum_{l=1}^{j-1} h^{(l)},p)
\eea
\end{lemma}

The Lemma follows from \Ref{De} and \Ref{D.op}.
\end{proof}

\subsection{Monodromy of solutions with respect to shifts of $z_j$ by $\tau$}
Let \fixedline
$u(z_1,\ldots,z_n,\la,\mu,\tau,p)$ be the function constructed for the tensor
product $V=V_{\La_1}\otimes\ldots\otimes V_{\La_n}$ in \Ref{sol}.
Let $D(\mu,z,p)$ be the $\End(V[0])$-valued function defined in \Ref{De}. Set
$$
\Psi(z,\la,\mu,\tau,p)\,=\,(1\otimes D(\mu,z,p))\,
u(z_1,\ldots,z_n,\la,\mu,\tau,p)\,.
$$
According to Theorem \ref{u-solu}, the $V[0]\otimes V[0]$-valued function
$\Psi$ is a solution of the qKZB equations with respect to the first factor,
\bean
\Psi(\ldots,z_k+p,\ldots)\,=\,(K_k(z,\la,\tau,p)\otimes 1)\,\,
\Psi(\ldots,z_k,\ldots),\qquad k=1,\ldots,n\,.
\notag
\eean

Let $B_j(z,\la,p)\in\End (V[0])$ be the linear operator introduced in \Ref{Bk}.
According to Theorem \ref{mon.shi}, for any $j=1,\ldots,n$, the function
$$
\Psi_j(z,\la,\mu,\tau,p)\,=\,
(B_j(z,\la, p)\otimes 1)\,\Psi(z_1,\dots,z_j+\tau,\dots,z_n,\la,\mu,\tau,p)
$$
is a new solution of the same equations.

The next Theorem describes a relation between the two solutions and can be
considered as a description of the monodromy of the hypergeometric solutions
constructed in Section \ref{sol.qKZB} with respect to shifts of variables
$z_j$ by $\tau$.

\begin{thm}\label{mondr.shifts}
\bean\label{m.s}
\kern-30pt&&
\Psi_j(z,\la,\mu,\tau,p)\,=\,\bigl(1\,\otimes\,f_j(z,p)
\\[2pt]
\kern-30pt&&
{}\times\,D(\ldots,z_j+\tau,\ldots,\mu,p)\,
K_j(z,\mu,p,\tau)\,D^{-1}(\ldots,z_j,\ldots,\mu, p)\bigr)\,
\Psi(z,\la,\mu,\tau,p)
\notag
\eean
where
$f_j(z,p)=e^{2\pi i\eta\sum_{l,\,l\neq j}(z_l-z_j)\Lambda_l\Lambda_j/p}$
and
$K_j(z,\mu,p,\tau)$ is the $j$-th operator of the qKZB equations with step
$\tau$ and modulus $p$.
\end{thm}

\noindent{\sl Remark.} According to Theorem \ref{tr.qkzb},
the operator
$$
\tilde K_j (z,\mu,p,\tau)\,=\,f_j(z,p)D(\ldots,z_j+\tau,\ldots,\mu,p)\,
K_j(z,\mu,p, \tau)\,D^{-1}(\ldots,z_j,\ldots,\mu, p)
$$
is $p$-periodic,
$\tilde K_j (\ldots,z_k+p,\ldots)=\tilde K_j (\ldots,z_k,\ldots)$.
Hence, formula \Ref{m.s} expresses solution $\Psi_j$ as a linear combination
of solutions $\Psi$ with $p$-periodic coefficients.

\vs
\begin{proof}
\bean
\kern-20pt&&
\Psi_j(z,\la,\mu,\tau,p)\,=\,
(B_j(z,\la,p)\otimes 1)\,\Psi(\ldots,z_j+\tau,\ldots)\,={}
\notag
\\
\kern-20pt&&
(f_j(z,p)\,D^{-1}(z,\la,\tau)D(\ldots,z_j+\tau,\ldots,\la,\tau)
\otimes D(\ldots,z_j+\tau,\ldots,\mu,p))\,u(\ldots,z_j+\tau,\ldots )\,={}
\notag
\\
\kern-20pt&&
(f_j(z,p)\,D^{-1}(z,\la,\tau)\,D(z,\la,\tau)\otimes
D(\ldots,z_j+\tau,\ldots,\mu,p)\,K_j(z,\mu,p, \tau))\,u(\ldots,z_j,\ldots)\,={}
\notag
\\
\kern-20pt&&
(1\otimes f_j(z,p)\,D(\ldots,z_j+\tau,\ldots,\mu,p)\,K_j(z,\mu,p,\tau)
\,D^{-1}(z,\mu,p))\,\Psi(\ldots,z_j,\ldots)\,.
\notag
\eean
\end{proof}

\noindent{\sl Remark.} Assume that $\La_1,\ldots,\La_n$ are natural numbers,
then the square of the transformation
\bean\label{tr.shift}
T_j:\Psi(z_1,\ldots,z_n,\la)\,\to\,
B_j(z,\la, p)\,\Psi(z_1,\dots,z_j+\tau,\dots,z_n,\la)
\notag
\eean
preserves the theta function properties \Ref{theta1}, \Ref{theta2}.
Namely, if a $V[0]$-valued function $\Psi$ satisfies \Ref{theta1},
\Ref{theta2}, then for any $j$, the $V[0]$-valued function $(T_j)^2\Psi$
satisfies \Ref{theta1}, \Ref{theta2}.

\subsection{Monodromy of solutions with respect to shifts of $z_j$ by $1$}
Let \fixedline
$u^\xi (z_1,\ldots,z_n,\la, \mu, \tau,p)$
be the function constructed for the tensor product
$V=V_{\La_1}\otimes\ldots\otimes V_{\La_n}$ in \Ref{sol}.
Let $D(\mu,z,p)$ be the $\End(V[0])$-valued function defined in \Ref{De}. Set
$$
\Psi^\xi(z,\la,\mu,\tau,p)\,=\,
(1\otimes D(\mu, z,p))\,u^\xi(z_1,\ldots,z_n,\la,\mu,\tau,p)\,.
$$
According to Theorem \ref{u-solu}, the $V[0]\otimes V[0]$-valued function
$\Psi^\xi$ is a solution of the qKZB equations with respect to the first
factor. According to Corollary \ref{shifts-1}, for any $j$,
the $V[0]\otimes V[0]$-valued function
$\Psi^\xi_j(z,\la,\mu,\tau,p)=\Psi(\ldots,z_j+1,\ldots,\la,\mu,\tau,p)$ is a
new solution of the same equations.

The next Proposition describes a relation between the two solutions and can be
considered as a description of the monodromy of the hypergeometric solutions
constructed in Section \ref{sol.qKZB} with respect to shifts of variables
$z_j$ by $1$.

For any $j$, introduce an entire function, $\xi_j$, of one variable by
$\xi_j(\la)=\xi(\la-2\eta \La_j)$.

\begin{proposition}\label{mondr.shifts-1}
For any $j$, we have
\bea
\Psi^\xi_j(z,\la,\mu,\tau,p) \, = \,(1\otimes (D_j(\mu))^{1/p})\,
\Psi^{\xi_j}(z,\la,\mu,\tau,p)
\eea
\end{proposition}
The Proposition follows from formula \Ref{coord}.

\vs
\noindent{\sl Remark.} Assume that $\La_1,\ldots,\La_n$ are natural numbers,
then the square of the transformation
\bean\label{tr.shift-1}
T_j:\Psi(z_1,\ldots,z_n,\la) \,\to\,\Psi(z_1,\dots,z_j+1,\dots,z_n,\la)
\notag
\eean
preserves the theta function properties \Ref{theta1}, \Ref{theta2}. Namely, if
a $V[0]$-valued function $\Psi (z,\la)$ satisfies \Ref{theta1}, \Ref{theta2},
then for any $j$, the $V[0]$-valued function $\Psi(\ldots,z_j+2,\ldots,\la)$
satisfies \Ref{theta1}, \Ref{theta2}.

Recall also that if $-p/4\eta=N$, $f$ is a scalar theta function of level $N$
and $\xi(\la)=f(\la/p)$, then $\Psi^\xi(z,\la,\mu,\tau,p)$ is a vector-valued
theta function of level $m+N$. Notice now that if $\La_1,\ldots,\La_n$ are
natural numbers, then $\xi(\la-4\eta\La_j)=f(\la/p+\La_j/N)$, and
$f(\la+\La_j/N)$ is a new scalar theta function of level $N$.

\end{document}